\documentclass[11pt]{amsart}

\usepackage{amsmath}
\usepackage{amssymb}
\usepackage{mathrsfs}
\usepackage[all]{xy}
\usepackage{dsfont}

\newtheorem{theorem}{Theorem}[section]
\newtheorem{claim}[theorem]{Claim}

\newtheorem{lemma}[theorem]{Lemma}

\newtheorem{proposition}[theorem]{Proposition}
\newtheorem{corollary}[theorem]{Corollary}

\theoremstyle{definition}
\newtheorem{definition}[theorem]{Definition}

\newtheorem{question}[theorem]{Question}

\theoremstyle{remark}
\newtheorem{remark}[theorem]{Remark}

\newtheorem{notation}[theorem]{Notation}

\newcommand{\dom}{\mathrm{dom}}
\newcommand{\gch}{\mathsf{GCH}}
\def\l{{\langle}}
\def\r{{\rangle}}
\def\s{\subseteq}
\DeclareMathOperator{\crit}{crit}
\DeclareMathOperator{\col}{\mathop{\mathrm{Col}}}
\DeclareMathOperator{\Add}{Add}
\DeclareMathOperator{\id}{id}
\DeclareMathOperator{\ran}{ran}
\def\br{\blacktriangleright}
\def\forces{\Vdash}
\newcommand{\len}{\mathrm{len}}
\newcommand{\Lim}{\mathrm{Lim}}
\newcommand{\supp}{\mathrm{supp}}
\newcommand\cat[1]{{}^\curvearrowright\langle #1\rangle}
\newcommand{\stem}{\mathrm{stem}}
\newcommand\diagonal{\bigtriangleup}

\newcommand{\Gal}[2]{\mathrm{Gal}(\mathrm{Cub}_{#1},{#1},{#2})}

\newcount\skewfactor
\def\mathunderaccent#1#2 {\let\theaccent#1\skewfactor#2
\mathpalette\putaccentunder}
\def\putaccentunder#1#2{\oalign{$#1#2$\crcr\hidewidth
\vbox to.2ex{\hbox{$#1\skew\skewfactor\theaccent{}$}\vss}\hidewidth}}
\def\name{\mathunderaccent\tilde-3 }

\newcommand{\lusim}[1]{\smash{\underset{\raisebox{1.2pt}[0cm][0cm]{$\sim$}}
{{#1}}}}
\def\smallbox#1{\leavevmode\thinspace\hbox{\vrule\vtop{\vbox
   {\hrule\kern1pt\hbox{\vphantom{\tt/}\thinspace{\tt#1}\thinspace}}
   \kern1pt\hrule}\vrule}\thinspace}

\newcommand{\cf}{{\rm cf}}

\title{Non-Galvin filters}

\author{Tom Benhamou}
\address[Benhamou]{Department of Mathematics, Statistics, and Computer Science, University of Illinois at Chicago, Chicago, IL 60607, USA}
\email{tomb@uic.edu}

\author{Shimon Garti}
\address[Garti]{Einstein Institute of Mathematics,
 The Hebrew University of Jerusalem,
 Jerusalem 91904, Israel}
\email{shimon.garty@mail.huji.ac.il}

\author{Moti Gitik}
\thanks{The research of the third author was partially supported by the by ISF grant 882/22.}
\address[Gitik]{School of Mathematical Sciences, Raymond and Beverly Sackler Faculty of Exact Science, Tel-Aviv University, Ramat Aviv 69978, Israel}
\email{gitik@post.tau.ac.il}
\author{Alejandro Poveda}
\address[Poveda]{Center of Mathematical Sciences and Applications, Harvard University, Cambridge MA, 02138, USA
 }
\email{alejandro@cmsa.fas.harvard.edu}


\begin{document}
\let\labeloriginal\label
\let\reforiginal\ref
\def\ref#1{\reforiginal{#1}}
\def\label#1{\labeloriginal{#1}}

\maketitle 
\begin{abstract}
    We address the question of consistency strength of certain filters and ultrafilters which fails to satisfy the Galvin property. We answer questions \cite[Questions 7.8,7.9]{TomMotiII}, \cite[Question 5]{NegGalSing} and improve theorem \cite[Theorem 2.3]{NegGalSing}. 
\end{abstract}
\section{introduction}
In this paper we continue the investigations on Galvin's property from \cite{TomMotiII, BenGarShe, NegGalSing,bgp}.  Let $\kappa$ be a regular uncountable cardinal and $\mathscr{F}$  a  $\kappa$-complete filter  over it. We shall write $\mathrm{Gal}(\mathscr{F},\kappa,\kappa^+)$ as a shorthand for the following statement:  every $\langle A_\alpha\mid \alpha<\kappa^+\rangle\s \mathscr{F}$ admits a subsequence $\langle A_{\alpha_\beta}\mid \beta<\kappa\rangle$ such that $\bigcap_{\beta<\kappa} A_{\alpha_\beta}\in\mathscr{F}$. If $\mathrm{Gal}(\mathscr{F},\kappa,\kappa^+)$ holds we shall say  that \emph{Galvin's property holds for $\mathscr{F}$} or, simply, that \emph{$\mathscr{F}$ is Galvin}. This terminology is coined in homage to F. Galvin's  discovery that if $\kappa^{<\kappa}=\kappa$ then   the \emph{club filter over $\kappa$}  ($\mathrm{Cub}_\kappa$) is Galvin \cite{MR0369081}. More generally, Galvin's proof shows that $\mathrm{Gal}(\mathscr{F},\kappa,\kappa^+)$ holds provided $\kappa^{<\kappa}=\kappa$ and $\mathscr{F}$ is normal.

The purpose of this paper is to present several constructions, both in the context of filters and ultrafilters, where Galvin's property fails. 
The first consistent example of a non-Galvin filter was provided by Abraham and Shelah \cite{MR830084}. In the said paper the authors exhibit a forcing poset producing a generic extension where $\mathrm{Gal}(\mathrm{Cub}_{\kappa^+},\kappa^+,{\kappa^{++}})$ fails for a regular cardinal $\kappa$. By virtue of Galvin's theorem,  $2^\kappa>\kappa^{+}$ in this latter model. An example of a ultrafilter $\mathscr{U}\s\mathcal{P}(\kappa)$ for which $\mathrm{Gal}(\mathscr{U},\kappa,\kappa^+)$ fails was given by Benhamou, Garti and Shelah \cite{BenGarShe}. Recently, in \cite{NegGalSing} it was shown how to make $\mathrm{Cub}_{\kappa^+}$ non-Galvin for all singular cardinal $\kappa$, simultaneously. 

\smallskip

The present manuscript is articulated in three blocks. In the first block (\S\ref{NoGalvinForUltra}) we analyze the failure of Galvin's property for ultrafilters that extend the club filter. This issue was first raised in \cite{TomMotiII} and subsequently answered in \cite{BenGarShe} under the existence of a supercompact cardinal. Shortly after this was improved in \cite{OnPrikryandCohen} using just a measurable cardinal. Here we modify the construction from \cite{BenGarShe} aiming to produce ultrafilters $\mathscr{U}$ concentrating on the set of singular cardinals, $\{\alpha<\kappa\mid \cf(\alpha)<\alpha\}$. This method is flexible-enough to generate $\kappa$-complete ultrafilters $\mathscr{U}$ such that $\mathrm{Cub}_\kappa\s \mathscr{U}$ and $\{\alpha<\kappa\mid \cf(\alpha)=\alpha\}\in \mathscr{U}$ (Theorem~\ref{thmregulars}). These are the sort of ultrafilters constructed in \cite{OnPrikryandCohen} using a completely different method. The advantage of the current  strategy in front the one of \cite{OnPrikryandCohen} is that the former, besides, adapts  to handle the singular case. Following up with this issue, in Theorem~\ref{ConSing} we give a lower bound for the consistency-strength of ``There is a $\kappa$-complete ultrafilter $\mathscr{U}\supseteq \mathrm{Cub}_\kappa$ with $\{\alpha<\kappa\mid \cf(\alpha)<\alpha\}\in \mathscr{U}$'' - this being $o(\kappa)\geq 2.$ Later, in Theorem~\ref{consingopt}, we show starting from $o(\kappa)=2$ (i.e., from optimal assumptions) that it is possible to force a $\kappa$-complete ultrafilter $\mathscr{U}$ as above for which $\mathrm{Gal}(\mathscr{U},\kappa,\kappa^+)$ fails. The  idea is to combine the Kurepa-tree-approach of \cite{BenGarShe} and Gitik's construction of a $\kappa$-complete ultrafilter concentrating on singular cardinals \cite{GitikFormerRegulars}.
  
\smallskip

In the second block of this paper (\S\ref{SectionNonGalvinFilters} and \S\ref{SectionFailurefortheClubfilter}) we focus on failures of Galvin's property for filters. We commence with \S\ref{SectionNonGalvinFilters} showing the consistency of \textsf{GCH} with every regular cardinal $\kappa$ carrying a $\kappa$-complete non-Galvin filter. In particular, the normality assumption in Galvin's theorem is nece\-ssary. The key idea here is that the existence of $\kappa$-independent families $\mathscr{F}\s \mathcal{P}(\kappa)$ (see page~\pageref{independentfamily}) yield such filters. It should be emphasized that we produce these configurations without bearing on any large-cardinal assumption. However, the disadvantage of this approach seems to be that  the  filters generated do not contain the club filter. We address this issue in \S\ref{SectionManyNonGalvinWithClubfilter} where we prove the consistency of every singular cardinal $\kappa$ carrying a $\kappa^+$-complete filter $\mathscr{F}$ such that $\mathrm{Cub}_{\kappa^+}\s\mathscr{F}$ and $\neg \mathrm{Gal}(\mathscr{F},\kappa^+,\kappa^{++})$. Unlike the previous approach, this latter consistency result uses large cardinals.

In \S\ref{SectionFailurefortheClubfilter}, we describe how to produce $\kappa$-complete ultrafilters $\mathscr{U}\s\mathcal{P}(\kappa)$ with $\mathrm{Cub}_{\kappa}\s \mathscr{U}$  and $\{\alpha<\kappa\mid \text{$\Gal{\alpha^+}{\alpha^{++}}$ fails}\}\in\mathscr{U}$. In particular, after Tree-Prikry-forcing with respect to $\mathscr{U}$ one gets a model where $\kappa$ is singular and there are cofinally many failures of Galvin's property below it.  This can be used to illustrate a sort of failure of \emph{compactness} at $\kappa$ relative to this property. In \S\ref{SectionTheModifiedPrikry} we take a slightly different approach and show how to produce a similar configuration for the first singular cardinal, $\aleph_\omega$. The idea here is to introduce a Prikry sequence on a measurable cardinal and, simultaneously, force with the poset of Abraham and Shelah from \cite{MR830084}.  The section ends indicating why Prikry-type forcings seem not useful to produce infinitely-many consecutive  failures of Galvin's property.

\smallskip

The third and last block (\S\ref{SectionConsistencyStrength}) deals with the consistency-strength  of the failure of Galvin's property at the successor of a singular cardinal. In Theorem~\ref{exactstrengthlocal} we show that $\neg \Gal{\aleph_{\omega+1}}{\aleph_{\omega+2}}$ is forceable starting with a cardinal $\kappa$ carrying a $(\kappa,\kappa^{++})$-extender. In particular, this pins down the consistency strength of this property to the optimal one; namely, $o(\kappa)=\kappa^{++}$. This answers a question from \cite[\S5]{NegGalSing}. In addition, we get a close-to-optimal upper bound for the consistency strength of $``\Gal{\kappa^+}{\kappa^{++}}$ fails for every singular cardinal $\kappa$''. Specifically, we show that this is forceable starting with a $(\kappa+3)$-strong cardinal. This improves \cite[Theorem~2.3]{NegGalSing}.



\subsection{Notation} Our notation is standard and mostly follows  \cite{NegGalSing}. We force in the Israel style where $p\leq q$ means that $q\Vdash p\in G$ (i.e., \emph{$q$ is stronger than $p$}). Often times, we also write $q\geq p$. The notation $\mathrm{Cub}_{\kappa}$ is reserved for the \textit{club filter} on $\kappa$. For two regular cardinals $\mu<\kappa$, 
$E^\kappa_{\mu}$ denotes the set of ordinals $\alpha<\kappa$ with $\cf(\alpha)=\mu$. The set of regulars below a cardinal $\kappa$ will be denoted by $\mathrm{Reg}_\kappa$. For ultrafilters $\mathscr{U}$ and $\mathscr{V}$ over $\kappa$ we write  $\mathscr{U}\leq_{\mathrm{RK}} \mathscr{V}$ whenever \emph{$\mathscr{U}$ is Rudin-Keisler below $\mathscr{V}$}; namely, if  there is a function $f\colon \kappa\rightarrow \kappa$ such that for every $X\s\kappa$, $X\in\mathscr{U}$ if and only if $f^{-1}[X]\in\mathscr{V}$.
\section{The failure of Galvin's property for ultrafilters}\label{NoGalvinForUltra}
\subsection{Non-Galvin ultrafilters}

Let $\kappa$ be a measurable cardinal.
In \cite{BenGarShe} $S$-slim Kurepa trees were used in order to force the existence of a $\kappa$-complete ultrafilter $\mathscr{U}$ over $\kappa$ which is not Galvin.
The idea is that Galvin's property yields some combinatorial property while $S$-slim Kurepa rules it out.
This helpful idea works, essentially, only if the tree is $S$-slim and $S\subseteq E^\kappa_\theta$ for some $\theta<\kappa$.
In such cases, the associated coloring has but $\theta$-many colors, and then the negative relation is meaningful.
If the number of colors is $\kappa$ then a negative relation is trivial, and the above argument breaks down.
Thus if $S$ is a stationary subset of $\mathrm{Reg}_\kappa$ and one wishes to force the failure of Galvin's property at a $\kappa$-complete ultrafilter which contains $\mathrm{Reg}_\kappa$ then some modification of the above idea is required. 

Such ultrafilters were already constructed in \cite{OnPrikryandCohen} from just a measurable cardinal, and our objective here  is merely to expand the method of \cite{BenGarShe} and to force $\neg{\rm Gal}(\mathscr{U},\kappa,\kappa^+)$ where $\kappa$ is measurable and $\mathscr{U}$ is a $\kappa$-complete ultrafilter over $\kappa$ which concentrates on $\mathrm{Reg}_\kappa$.
We shall do it by modifying the forcing of \cite{BenGarShe}, but for our argument we need, first of all, a simple observation regarding the Rudin-Keisler order.
\begin{lemma}
\label{lemrk} Assume that:
\begin{enumerate}
\item [$(\aleph)$] $\mathscr{U},\mathscr{V}$ are ultrafilters over $\kappa$.
\item [$(\beth)$] $\mathscr{U}\leq_{\rm RK}\mathscr{V}$.
\item [$(\gimel)$] ${\rm Gal}(\mathscr{V},\kappa,\kappa^+)$.
\end{enumerate}
Then ${\rm Gal}(\mathscr{U},\kappa,\kappa^+)$.
\end{lemma}
\begin{proof}
Fix $\pi:\kappa\rightarrow\kappa$ witnessing the assumption $\mathscr{U}\leq_{\rm RK}\mathscr{V}$.
Suppose that $\l C_i\mid i<\kappa^+\}\subseteq\mathscr{U}$.
By definition, $\l\pi^{-1}[C_i]\mid i<\kappa^+\}\subseteq\mathscr{V}$, so one can find $I\in[\kappa^+]^\kappa$ and $B\in\mathscr{V}$ so that $B\subseteq\bigcap_{i\in I}\pi^{-1}[C_i]$.
Let $A=\pi{''}{B}$.
Notice that $A\subseteq{C_i}$ for every $i\in{I}$, thus ${\rm Gal}(\mathscr{U},\kappa,\kappa^+)$ is established.
\end{proof}

From the above lemma we infer that if one forces $\neg{\rm Gal}(\mathscr{U},\kappa,\kappa^+)$ and $\mathscr{U}\leq_{\rm RK}\mathscr{V}$ then $\neg{\rm Gal}(\mathscr{V},\kappa,\kappa^+)$.
Our strategy will be to force this situation where $\mathscr{V}$ concentrates on $\mathrm{Reg}_\kappa$.
This will be done by adding one feature to the forcing construction of \cite{BenGarShe}.
We recall the definition of the forcing $\mathbb{K}(S)$, where $S$ is a stationary subset of $\kappa$.

This forcing notion consists of two components.
The first one, $\mathbb{Q}(S)$, adds a stationary subset of $S$.
The second adds an $S$-slim Kurepa tree.
Thus $\mathbb{Q}(S)=\{x\subseteq{S}:|x|<\kappa\}$ and the order is end-extension.
The forcing notion $\mathbb{K}(S)$ consists of triples $(x,t,f)$ where $x\in\mathbb{Q}(S)$, $t$ is a normal tree of height $\beta+1$ for some $\beta\in\kappa$, $\sup(x)\geq\beta+1$ and $|\mathcal{L}_\alpha(t)|\leq|\alpha|$ whenever $\alpha\in x\cap\beta+1$.
Finally, $f:\kappa^+\rightarrow\mathcal{L}_\beta(t)$ is a partial function with $|f|\leq|\beta|$.
If $(x,t,f),(y,s,g)\in\mathbb{K}(S)$ then $(x,t,f)\leq_{\mathbb{K}(S)}(y,s,g)$ iff $x\subseteq_{\rm end}y, s\upharpoonright(\beta+1)=t, {\rm dom}(f)\subseteq{\rm dom}(g)$ and $f(\alpha)\leq_sg(\alpha)$ for every $\alpha\in{\rm dom}(f)$.
If $G\subseteq\mathbb{K}(S)$ if generic then $\mathscr{T}_G=\bigcup\{t:\exists x,f, (x,t,f)\in{G}\}$ is the desired slim tree.

\begin{theorem}
\label{thmregulars} Let $\kappa$ be supercompact.
Then one can force $\neg{\rm Gal}(\mathscr{V},\kappa,\kappa^+)$ where $\mathscr{V}$ is a $\kappa$-complete ultrafilter over $\kappa$ such that $\mathrm{Cub}_\kappa\subseteq \mathscr{V}$ and $\mathrm{Reg}_\kappa\in\mathscr{V}$.
\end{theorem}
\begin{proof}
Assume that $\kappa$ is supercompact and
let $h:\kappa\rightarrow{V_\kappa}$ be a Laver-diamond function. Let $S=E^\kappa_\omega$,
we define a two-step iteration $\mathbb{S}=\mathbb{P}\ast\name{\mathbb{R}}$ as follows.
The first component $\mathbb{P}$ is an Easton support iteration $\langle\mathbb{P}_\alpha,\name{\mathbb{Q}}_\beta: \alpha\leq\kappa,\beta<\kappa\rangle$, where $\name{\mathbb{Q}}_\beta$ is trivial unless $\beta$ is strongly inaccessible and $h(\beta)$ is a $\mathbb{P}_\beta$-name of the forcing notion $\mathbb{K}(S^\beta_\omega)\times Add(\beta,1)$, in which case we let $\name{\mathbb{Q}}_\beta=h(\beta)$.
The second component $\name{\mathbb{R}}$ is (a $\mathbb{P}$-name of) the forcing notion $\mathbb{K}(S)\times Add(\kappa,1)$.

Let $G\subseteq\mathbb{S}$ be $V$-generic, so $G$ factors into $G_{\mathbb{P}}\ast G_{\mathbb{R}}$ in a natural way.
From \cite{BenGarShe} we know that in $V[G]$ there is a stationary $S\subseteq E^\kappa_\omega$ and an $S$-slim Kurepa tree.
Moreover, there is a $\kappa$-complete ultrafilter $\mathscr{U}$ over $\kappa$ which extends $\mathrm{Cub}_\kappa\cup\{S\}$.
Let us briefly describe $\mathscr{U}$, and build another $\kappa$-complete ultrafilter $\mathscr{V}$ so that $\mathscr{U}\leq_{\rm RK}\mathscr{V}$ and $\mathrm{Reg}_\kappa\in\mathscr{V}$.

Choose $\lambda>2^\kappa$ and a supercompact elementary embedding $\jmath:V\rightarrow{M}$ such that ${\rm crit}(\jmath)=\kappa, \jmath(\kappa)=\lambda$ and ${}^{2^\kappa}M\subseteq{M}$. 
We also require that $\jmath(h)(\kappa)=\name{\mathbb{R}}$.
Let $\mathbb{P}'=\jmath(\mathbb{P})$, so $\mathbb{P}'= \langle\mathbb{P}'_\alpha,\name{\mathbb{Q}}'_\beta: \alpha\leq\jmath(\kappa),\beta<\jmath(\kappa)\rangle$.
Up to $\kappa$ we know that $\mathbb{P}'$ coincides with $\mathbb{P}$, since $\kappa={\rm crit}(\jmath)$.
We also know that $\mathbb{P}'_{\kappa+1}=\mathbb{P}\ast\name{\mathbb{R}}$ since $\jmath(h)(\kappa)=\name{\mathbb{R}}$ and by the definition of our iteration.

In particular, one can form the generic extension $M[G]$ in $V[G]$.
Observe that the rest of the iteration, that is, $\mathbb{P}'_{(\kappa+1,\jmath(\kappa))}$, is $\theta$-closed where $\theta$ is the first $M[G]$-inaccessible above $\kappa$.
In particular, it is $(2^\kappa)^+$-closed, where $(2^\kappa)^+$ is computed in $V$.

Working in $V$, let $\mathcal{C}=\{\name{C}:\name{C}$ 
is a nice $\mathbb{K}(S)$-name for a club of $\kappa\}$, so $|\mathcal{C}|=2^\kappa$.
Similarly, let $\mathcal{A}$ be $\{\name{A}:\name{A}$ 
is a nice $\mathbb{K}(S)$-name for a subset of $\kappa\}$ and then $|\mathcal{A}|=2^\kappa$.
Since ${}^{2^\kappa}M\subseteq{M}$ we see that both $\{\jmath(\name{C}):\name{C}\in\mathcal{C}\}$ and $\{\jmath(\name{A}):\name{A}\in\mathcal{A}\}$ are elements of $M$.

As a first step towards the construction of $\mathscr{U}$ we claim that there exist an ordinal $\delta$ and a condition $p\in\mathbb{P}'_{(\kappa+1,\jmath(\kappa))}$ such that $p$ forces in ${\mathbb{P}'_{(\kappa+1,\jmath(\kappa))}}$ that $\delta\in\bigcap\{\jmath(\name{C}): \name{C}\in\mathcal{C}\}\cap\jmath(\name{S})$.
To see this, recall that each $\jmath(\name{C})$ is a club of $\jmath(\kappa)$, and $\jmath(\kappa)=\lambda>2^\kappa$.
Thus, $\mathcal{C}$ is a collection of $2^\kappa$ clubs of $\jmath(\kappa)$ and hence it is forced by the empty condition that $\bigcap\{\jmath(\name{C}): \name{C}\in\mathcal{C}\}$ is a club of $\jmath(\kappa)$.
In addition, $\jmath(\name{S})$ is a stationary subset of $\jmath(\kappa)$, so one can find $\delta$ and $p$ such that $p\Vdash_{\mathbb{P}'_{(\kappa+1,\jmath(\kappa))}} \delta\in\bigcap\{\jmath(\name{C}): \name{C}\in\mathcal{C}\}\cap\jmath(\name{S})$.

By the closure of $\mathbb{P}'_{(\kappa+1,\jmath(\kappa))}$ there is $q\geq{p}$ such that $q$ decides the statement $\delta\in\jmath(\name{A})$ for every $\name{A}\in\mathcal{A}$.
Indeed, enumerate $\mathcal{A}$ by $\{\name{A}_j:j\in 2^\kappa\}$ and create an increasing sequence of conditions $\l q_j\mid j< 2^\kappa\r$ such that $q_j\parallel\delta\in\name{A}_j$ and $p\leq q_0$.
At the end, let $q$ be an upper bound of every $q_j$.

Now in $V[G]$ we can define $\mathscr{U}$ as the set $\{(\name{A})_G:\name{A}\in\mathcal{A}\wedge q\Vdash\delta\in\jmath(\name{A})\}$.
One can verify that $\mathscr{U}$ is a $\kappa$-complete ultrafilter over $\kappa$ (using the fact that $\kappa={\rm crit}(\jmath)$ and the elementarity of $\jmath$).
Moreover, $\mathrm{Cub}_\kappa\cup\{S\}\subseteq\mathscr{U}$ by the choice of $\delta$.
Our goal, therefore, is to construct $\mathscr{V}$.

Recall that $\name{\mathbb{R}}$ is $\mathbb{K}(S)\ast{Add(\kappa,1)}$, so let $f:\kappa\rightarrow\kappa$ be the Cohen part as interpreted by $G_{\mathbb{R}}$.
We claim that there are a condition $r\geq{q}$ and an $M[G]$-inaccessible $\rho\in\jmath(\kappa)$ such that:
\[r\Vdash_{\mathbb{P}'_{(\kappa+1,\jmath(\kappa))}} \rho\in\bigcap\{\jmath(\name{C}): \name{C}\in\mathcal{C}\} \wedge \jmath(\name{f})(\rho)=\delta\]
The claim is justified by the general fact that if $\mu$ is Mahlo in some model $W$ of \textsf{ZFC} and $g:\mu\rightarrow\mu$ is $W$-generic for the forcing notion $Add(\mu,1)$ then for every $\gamma<\mu$, the set $\{\alpha<\mu:\alpha$ is inaccessible and $g(\alpha)=\gamma\}$ is a stationary subset of $\mu$ in $W[g]$.

Our $\jmath(\kappa)$ is certainly Mahlo in $M[G]$, thus we may apply the general fact to $M[G]$ and $\jmath(f)$, and deduce that for every $\gamma<\jmath(\kappa)$ the condition $q$ forces that $S_\gamma=\{\alpha<\jmath(\kappa):\alpha$ is inaccessible and $\jmath(f)(\alpha)=\gamma\}$ is a stationary subset of $\jmath(\kappa)$.
Taking $\gamma$ as our $\delta$ we see that $S_\delta$ is stationary in $\jmath(\kappa)$, so there are $\rho\in{S_\delta}$ and $r\geq{q}$ such that $r\Vdash\jmath\name{f}(\rho)=\delta \wedge \rho\in\bigcap\{\jmath(\name{C}): \name{C}\in\mathcal{C}\}$. Again, by closure we can assume that such a condition already determines all the statements $\rho\in j(\name{A})$ for every $\name{A}\in \mathcal{A}$.
This choice enables us to define, in $V[G]$, the following set:
\[\mathscr{V}=\{(\name{A})_G\mid\name{A}\in\mathcal{A}, r\Vdash_{\mathbb{P}'_{(\kappa+1,\jmath(\kappa))}} \rho\in\jmath(\name{A})\}\]
It is routine to check that $\mathscr{V}$ is a $\kappa$-complete ultrafilter over $\kappa$ in $V[G]$.
Let us show that $\mathscr{U}\leq_{\rm RK}\mathscr{V}$ as witnessed by $f$.

To prove this fact, fix $A\in\mathscr{U}$, so $A=\name{A}_G$ for some $\name{A}\in\mathcal{A}$.
Let us show that $f^{-1}[A]\in\mathscr{V}$ (the opposite direction from $\mathscr{V}$ to $\mathscr{U}$ is similar).
By definition, $A\in\mathscr{U}$ implies that $q\Vdash\delta\in\jmath(\name{A})$.
Therefore, $q$ forces $\jmath(\name{f})(\rho)\in\jmath(\name{A})$ by the choice of $\rho$.
This means that $q$ forces $\rho\in\jmath(\name{f}^{-1}[\name{A}])$.
By the definition of $\mathscr{V}$ we conclude that $f^{-1}[\name{A}_G]\in\mathscr{V}$, thus $f^{-1}[A]\in\mathscr{V}$ as required.

Since $\rho$ is inaccessible in $M$, the set $\mathrm{Reg}_\kappa$ belongs to $\mathscr{V}$.
From Lemma \ref{lemrk} we know that ${\rm Gal}(\mathscr{V},\kappa,\kappa^+)$ holds true in $V[G]$, so we are done.
\end{proof}

Let us indicate that stronger properties can be forced upon $\mathscr{V}$, due to the choice of $\rho$.
Thus, since $\jmath(\kappa)$ is supercompact in $M$ one can choose $\rho$ to be measurable and then $\mathscr{V}$ concentrates on measurable cardinals.

\subsection{A non Galvin ultrafilter concentrating on singulars}
In this section we produce  a $\kappa$-complete non-Galvin ultrafilter that  concentrates on singulars and extends the club filter $\mathrm{Cub}_\kappa$. We will accomplish the construction starting from optimal large-cardinal assumptions, hence improving the main result of \cite{BenGarShe}.  The readers familiar with \cite{OnPrikryandCohen} will note that the present context differs from the former in that all the ultrafilters considered in \cite{OnPrikryandCohen} concentrated on the set of regular cardinals.

As the forthcoming theorem shows $o(\kappa)=2$ is the minimal large-cardinal assumption for the existence of a $\kappa$-complete ultrafilter $\mathscr{U}\supseteq \mathrm{Cub}_\kappa$ concentrating on the set of singular cardinals. The argument is due (basically) to W. Mitchell but we add the proof for the reader's convenience.
\begin{theorem}\label{ConSing}
 Suppose there is a $\kappa$-complete ultrafilter $U$ over $\kappa$ such that $\{\alpha<\kappa\mid \alpha\text{ is singular }\}\in U$ and $\mathrm{Cub}_\kappa\subseteq U$ or alternatively, that $[id]_U$ is a generator of $j_U$ and $[\id]_U$ is singular.  Then either there is an inner model with a Woodin cardinal or in the core model $\mathcal{K}$,  $o^{\mathcal{K}}(\kappa)\geq 2$.
 
\end{theorem}
\begin{proof}
 Suppose that there is no inner model with a Woodin cardinal and let $\mathcal{K}$ be the Jensen-Steel core model \cite{jensen_steel_2013}. Consider $j_U:V\rightarrow M_U$ the ultrapower embedding. By Schindler \cite{schindler_2006}  $j_U\restriction \mathcal{K}:\mathcal{K}\rightarrow\mathcal{K}^{M_U}$ is an iterated ultrapower of $\mathcal{K}$ by its measures $\l i_{\alpha,\beta}\mid \alpha\leq\theta\r$. Suppose that the iteration is normal and let $\l \kappa_i\mid i\leq\lambda\r$ be the increasing enumeration of $\{i_{0,\alpha}(\kappa)\mid \alpha\leq\theta\}$. Since $\mathrm{Cub}_\kappa\subseteq U$, there is a  $\delta<\lambda$ such that $\kappa_\delta=[id]_U$. Just otherwise $\kappa_\delta<[id]_U<\kappa_{\delta+1}$ and by \cite[Claim 44]{YairTomMoti}, there would be $f:\kappa\rightarrow\kappa$ such that $j_U(f)(\kappa_\delta)\geq [id]_U$. But then $C_f:=\{\alpha<\kappa\mid f''\alpha\subseteq\alpha\}$ would be a club at $\kappa$ and $[id]_U\notin j_U(C_f)$. Since $[id]_U$ is singular in $M_U$, it follows that $[id]_U=\kappa_\delta$ for some limit $\delta$ since by \cite[Lemma 46]{YairTomMoti} each successor element of the sequence of the form $\kappa_{i+1}$ is regular in $M_U$.
 
 To deduce that $o^{\mathcal{K}}(\kappa)\geq 2$, suppose otherwise that $o^{\mathcal{K}}(\kappa)=1$, and denote by $W$ the only measure on $\kappa$ in $\mathcal{K}$. Since $M_U$ is closed under $\omega$-sequences $\l\kappa_n\mid n<\omega\r\in M_U$. Now all the $\kappa_n$'s are critical points of the iteration, namely for some $\alpha_n$, $crit(i_{\alpha_n,\theta})=\kappa_n=i_{0,\alpha_n}(\kappa)$ (see for example \cite[Corollary 43]{YairTomMoti}) and $i_{\alpha_n,\alpha_{n+1}}$ is the ulrtapower embedding by $i_{0,\alpha_n}(W)$. Let $\alpha_\omega=\sup_{n<\omega}\alpha_n$. Note that $i_{0,\alpha_\omega}(\kappa)<i_\theta(\kappa)=j_U(\kappa)$, otherwise $cf^{M_U}(j_U(\kappa))=\omega$ which contradicts the elementarity of $j_U$. Hence by the normality of the iteration, $crit(i_{\alpha_\omega,\alpha_\omega+1})=i_{0,\alpha_\omega}(\kappa)$ and since $o(i_{0,\alpha_\omega}(\kappa))=1$ it follows that $i_{\alpha_\omega,\alpha_\omega+1}:\mathcal{K}_{\alpha_\omega}\rightarrow\mathcal{K}_{\alpha_\omega+1}$ is the ultrapower by $i_{0,\alpha_\omega}(W)$. In particular, $i_{0,\alpha_\omega}(W)\notin \mathcal{K}_{\alpha_\omega+1}$. Working in $M_U$, using the sequence of $\kappa_n$'s which forms a Prikry sequence for $i_{0,\alpha_\omega}(W)$, we can reconstruct $i_{0,\alpha_\omega}(W)\in M_U$. Thus $i_{0,\alpha_\omega}(W)\in \mathcal{K}^{M_U}$ as any $\mathcal{K}^{M_U}$-measure in $M_U$ already belongs to $\mathcal{K}^{M_U}$. However, $i_{\alpha_\omega+1,\theta}:\mathcal{K}_{\alpha_\omega+1}\rightarrow\mathcal{K}^{M_U}$, and by normality of the iteration, $crit(i_{\alpha_\omega+1,\theta})$ is much above $i_{0,\alpha_\omega}(\kappa)$ which ensures that $i_{0,\alpha_\omega}(W)\in \mathcal{K}_{\alpha_\omega+1}$, contradiction.
\end{proof}
Let us prove that under the minimal assumption it is consistent to produce a witness for the negation of the Galvin property.
\begin{theorem}\label{consingopt}
 Assume $\mathsf{GCH}$ and suppose $o(\kappa)\geq 2$ then it is consistent that there is a $\kappa$-complete ultrafilter $W$ such that  $\mathrm{Cub}_\kappa\subseteq W$ and $\{\alpha<\kappa\mid \alpha\text{ is singular }\}\in W$ which fails to satisfy the Galvin property.
\end{theorem}
The main Lemma is the following:
\begin{lemma}\label{Preperation}
Suppose that $U_0\triangleleft U_1$ are normal measures over $\kappa$. Then there is a forcing extension $V[G]$ such that in $V[G]$ $U_0,U_1$ extend to $U_0^*,U_1^*$ respectively, $U_0^*$ is normal,  $U_0^{*}\leq_{RK}U^*_1$, $U_1^*$ concentrates on $E^\kappa_\omega$ and the $\omega$-iteration by $U^*_0$ denoted by  $\l j_{n,m},M_n\mid i<\omega\r$ satisfies the following:
\begin{enumerate}
    \item $\kappa=crit(j_{0,1})$, and $crit(j_{n,n+1})=j_{0,n}(\kappa)$. 
    \item $\l j_{0,n}(\kappa)\mid n<\omega\r$ is unbounded in $[id]_{U^*_1}$.
    \item There are factor maps for the embedding ultrapower $j_{U^*_1}$, $k_n: M_n\rightarrow M_{U^*_1}$ such that $j_{U^*_1}=k_n*j_{0,n}$, $k_n=k_m\circ j_{n,m}$ and $crit(k_n)=j_{0,n}(\kappa)$.
\end{enumerate}
\end{lemma}
Let us first conclude Theorem \ref{consingopt} from this lemma.

\textit{Proof of Theorem \ref{consingopt}.} Over $N:=V[G]$ we force the same forcing as in \cite{BenGarShe}, iterating with Easton support the forcing for adding $E^\alpha_\omega$-slim Kurepa trees for each $\alpha\leq\kappa$ denoted by $\mathbb{P}*\mathbb{R}$. Let $H_\kappa*h$ be an $N$-generic filter for 
$\mathbb{P}*\mathbb{R}$. First we extend $j_{0,1}:N\rightarrow M_{1}$, note that by Easton support, $$j_{0,1}(\mathbb{P}*\mathbb{R})=\mathbb{P}*\mathbb{R}*\mathbb{P}_{(\kappa+1,j_{0,1}(\kappa))}*j_{0,1}(\mathbb{R})$$ and $j_{0,1}''H_\kappa=H_\kappa$. Let us define in $N$ an $M_{1}$-generic filter for $j(\mathbb{P}*\mathbb{R})$ by first 
taking $H_{\kappa}*h$. Note that the forcing $\mathbb{P}_{(\kappa+1,j_{0,1}(\kappa))}$ starts above $\kappa^+$ and by $\textsf{GCH}$, there are only $\kappa^+$-many dense subsets to meet. By standard arguments, exploiting the fact that $M_{1}$ is the ultrapower by a $\kappa$-complete measure, hence closed under $\kappa$-sequences of 
$V$, we construct an $M_{1}[H_{\kappa}*h]$-generic filter $T$ for  $\mathbb{P}_{(\kappa+1,j_{0,1}(\kappa))}$. 

Notice that the model $M_{1}[H_{\kappa}*h*T]$ is closed only under $\kappa$-sequences from $N$ and since $|h|=\kappa^+$, we cannot guarantee that $j''h\in M_{1}[H_{\kappa}*h*T]$. Instead, we start by constructing any $M_{1}[H_{\kappa}*h*T]$-generic $t'$ for $j_{1}(\mathbb{R})$  starting above the condition $\l f_{\kappa},T_{\kappa}\r$ where $T_{\kappa}$ is a tree of height $\kappa+1$, $T_{\kappa}\restriction \kappa$ is constructed from $h$ and $Lev_{\kappa}(T_{\kappa})=\{b_\kappa(\alpha)\mid\alpha<\kappa^+\}$ where  $b_\kappa(\alpha)$ are the branches derived from $h$. It is crucial here that we can take $\kappa^+$-many elements in $Lev_{\kappa}(T_{\kappa})$ since $\kappa\notin j_{0,1}(E^{\kappa}_{\omega})=E^{\kappa_1}_\omega$. The function $f_\kappa:j_{0,1}(\kappa)^+\rightarrow Lev_{\kappa}(T_\kappa)$ is the partial function with $\dom(f_{\kappa})=\kappa$ and $f_{\kappa}(\alpha)=b_\kappa(\alpha)$ for every $\alpha<\kappa$.  Now in $N$ we define an additional filter $t$ which is obtained from $t'$ by changing the values of each branch of the form $b_{j_{0,1}(\kappa)}(j_{0,1}(\alpha))$ where $\alpha<\kappa^+$ so that $$b_{j_{0,1}(\kappa)}(j_{0,1}(\alpha))\restriction\kappa=b_\kappa(\alpha).$$ Formally, for every pair $\l f_\kappa,T_\kappa\r\leq\l g,S\r$ we define $\l g^*,S\r$ where for every $j_{0,1}(\alpha)\in \dom(g)\cap j_{0,1}''\kappa^+$ we let $g^*(j_{0,1}(\alpha))\restriction \kappa=b_{\kappa}(\alpha)$. Note that $g^*\in M_{1}[H_{\kappa}*h*T]$, since $\dom(g)$ is bounded in $j_{0,1}(\kappa)^+$, $|\dom(g)\cap j''\kappa^+|\leq\kappa$ and  $M_{1}[H_{\kappa}*h*T]$ is closed under $\kappa$-sequences. Define $t=\{\l b^*,B\r\mid \l b,B\r\in t'\}$, then $t\subseteq j_{0,1}(\mathbb{R})$ and it is $M_{1}[H_{\kappa}*h*T]$-generic. Indeed for every dense open set $D\in M_{1}[H_{\kappa}*h*T]$ we can find $D^*\in M[H_{\kappa}*h*T]$ still dense open such that for every $\l g,S\r\in D^*$ and any $\kappa$-many changes of $g$, $\l g^*,S\r\in D$. Let us denote $$H_{\kappa}=H_0, \ h=h_0, \ G_0=H_{0}*h_0\text{ and}$$  $$H_1=H_{\kappa}*h*T, \ h_1=t, \ G_1=H_1*h_1.$$
In the same fashion, keep defining inductively the generic filters $G_n=H_n*h_n$ such that $G_{n+1}\restriction \kappa_n+1=G_n$. 
So far we have extended the embeddings $j^*_{n,m}$ to form a directed system $\l j^*_{n,m},M_n[G_n]\mid n\leq m<\omega\r$. Denote the direct limit of the models $M_\omega=\underrightarrow{Lim}M_n$, the direct limit embeddings $j_{n,\omega}:M_n\rightarrow M_\omega$, and the direct limit of the factor maps $k_\omega:M_\omega\rightarrow M_{U^*_1}$, which is defined by the equalities $k_\omega\circ j_{n,\omega}=k_n$. Finally, denote  $\underrightarrow{Lim}M_n[G_n]=M_{\omega}[G_{\omega}]$ and the direct limit embeddings $j^*_{n,\omega}:M_n[G_n]\rightarrow M_\omega[G_\omega]$. Note that $$G_{\omega}=H_\omega*h_\omega=j_{n,\omega}^*(H_n*h_n).$$ 

 To extend $j_{U^*_1}$, we start the construction up to $[id]_{U^*_1}$ not including $[id]_{U^*_1}$. Recall that by the lemma, $crit(k_n)=j_{0,n}(\kappa)$ and $\l j_{0,n}(\kappa)\mid n<\omega\r$ is unbounded in $[id]_{U^*_1}$. It follows that $k_n''H_n=H_n$ and $H_\omega=\cup_{n<\omega}H_n$. Since $crit(k_\omega)=\sup j_{0,n}(\omega)=[id]_{U^*_1}$, it follows that $H_{\omega}$ is $M_{U^*_1}$-generic for $\mathbb{P}_{[id]_{U^*_1}}$. 
 
 Next, at $[id]_{U^*_1}$ which is of cofinality $\omega$ in $M_{U^*_1}$ the iteration is defined to be trivial. Finally, above $[id]_{U^*_1}$ we construct the generic in a similar fashion to what we did for $j_{U^*_0}$ where in the construction of the tree at $j_{U^*_1}(\mathbb{R})$ we change the values of the generic at $j_{U^*_1}(\kappa)$ with respect to the point-wise image by $k_\omega$ of $h_{\omega}$.
 
We obtain an $M_{U^*_1}$-generic filter $G^*\in N[G_0]$ and extend $j_{U^*_1}$ to $j^*_{U^*_1}:N[G_0]\rightarrow M_{U^*_1}[G^*]$. Note that by the construction of the generic filters $G_n$, also $k_n$ extends to commutative factor maps $k^*_n:M_n[G_n]\rightarrow M_{U^*_1}[G^*]$. Hence by the universal property of direct limits, there is $k^*_\omega:M_{\omega}[G_\omega]\rightarrow M_{U^*_1}[G^*]$ such that $k^*_\omega\circ j^*_{n,\omega}=k^*_n$. Since each $k^*_n$ extends $k_n$ we see that $k^*_\omega$ extends $k_\omega$.

In $N[G_0]$, derive the ultrafilter $W=\{X\subseteq\kappa\mid [id]_{U^*_1}\in j^*_{U^*_1}(X)\}$.

\begin{proposition}\label{thisprop}
$\mathrm{Cub}_\kappa\subseteq W$, $\{\alpha\mid cf(\alpha)=\omega\}\in W$ and $W$ fails to satisfy the Galvin property.
\end{proposition}
\textit{
Proof of Proposition \ref{thisprop}.}
 Since $M_{U^*_1}[G^*]\models cf([id]_{U^*_1})=\omega$, 
 $E^\kappa_\omega\in W$. Let $C\in \mathrm{Cub}_\kappa$, we would like to prove that $[id]_{U^*_1}\in j^*_{U^*_1}(C)$. By elementarity, $j^*_{U^*_1}(C)$ is closed and since $\l j_{0,n}(\kappa)\mid n<\omega\r$ is unbounded in $[id]_{U^*_1}$, it suffices to prove the for every $n<\omega$, $j_{0,n}(\kappa)\in j^*_{U^*_1}(C)$. Let $n<\omega$, since $crit(k^*_n)=crit(k_n)=j_{0,n}(\kappa)$, and $k^*_n(j^*_{0,n}(C))=j^*_{U^*_1}(C)$, $j^*_{U^*_1}(C)\cap j_{0,n}(\kappa)=j^*_{0,n}(C)$. By elementarity, $j^*_{0,n}(C)$ is unbounded in $j_{0,n}(\kappa)$, and since $j^*_{U^*_1}(C)$ is closed, $j_{0,n}(\kappa)\in j^*_{U^*_1}(C)$.
 Recall that in $N[G_0]$ we have an $E^\kappa_\omega$-slim Kurepa tree, so by \cite{BenGarShe}, $W$ cannot have the Galvin property.
 $\qed_{\text{Prop.}\ref{thisprop}}$
$\qed_{\text{Thm.}\ref{consingopt}}$

\textit{Proof of Lemma \ref{Preperation}}:
 Let $U_0\triangleleft U_1$ and $[\xi\mapsto U_0(\xi)]_{U_1}=U_0$. Let $$E=\{\alpha<\kappa\mid \alpha\text{ is measurable with measure }U_0(\alpha)\}\in U_1\setminus U_0.$$ Suppose also that $\forall \alpha\in E,E\cap\alpha\notin U_0(\alpha)$. Over the ground model $V$, we first force with $\mathbb{P}_\kappa$ which is an Easton support iteration of Prikry forcing as defined in \cite{584c9093529345ed9ca7a01a0e0c9d9a} adding for each $\alpha\in E$ a Prikry sequence $\{\alpha_n\mid n<\omega\}$. Namely:
 \begin{definition}
Let $\overline{E}$ be the closure of the set $E\cup\{\alpha+1\mid \alpha\in E\}\cup\{\kappa\}$. For $\alpha\in \overline{E}$, we inductively define $\mathbb{P}_\alpha$,  the conditions of $\mathbb{P}_\alpha$ are functions $p=\{\l t^p_\gamma, A^p_\gamma\r\mid \gamma\in \dom(p)\}$ such that:
 \begin{enumerate}
     \item $\dom(p)\subseteq \overline{E}\cap\alpha$.
     \item $p$ has Easton support, namely,  for every inaccessible $\beta\leq\alpha$, $\dom(p)\cap\beta$ is bounded in $\beta$.
     \item For every $\gamma\in \dom(p)$, $p\restriction\gamma:=\{\l t^p_\beta, A^p_\beta\r\mid \beta\in \dom(p)\cap\gamma\}\in \mathbb{P}_\gamma$ 
     and $p\restriction\gamma\Vdash_{\mathbb{P}_\gamma}\l t^p_\gamma,A^p_\gamma\r\in \mathcal{P}(\lusim{U}^*_0(\alpha))$, where $\lusim{U}^*_0(\alpha)$ is a normal ultrafilter over $\alpha$ defined in Definition \ref{definition of ultrafilter extension}. 
 \end{enumerate}
 The order is $p\leq q$ if and only if:
\begin{enumerate}
    \item $\dom(p)\subseteq\dom(q)$.
    \item For every $\gamma\in \dom(p)$, $q\restriction \gamma\Vdash_{\mathbb{P}_\gamma} \l t^p_\gamma,A^p_\gamma\r\leq  \l t^q_\gamma,A^q_\gamma\r$.
    \item There is a finite set $b$ such that for every $\gamma\in \dom(p)\setminus b$, $q\restriction \gamma\Vdash_{\mathbb{P}_\gamma} t^p_\gamma= t^q_\gamma$.
\end{enumerate} 
Moreover in clause $3.$ if $b=\emptyset$ then we say that $q$ is a direct extension of $p$ and denote it by $p\leq^*q$.

 \end{definition}
The following lemmas can be found in \cite{584c9093529345ed9ca7a01a0e0c9d9a}:
\begin{lemma}
Let $\l p_\beta\mid \beta<\gamma<\alpha\r$ be a sequence of conditions in $\mathbb{P}_\alpha$ such that for every $\beta_1\leq\beta_2$, $p_{\beta_1}\restriction \gamma+1=p_{\beta_2}\restriction\gamma+1$ and $p_{\beta_1}\leq^* p_{\beta_2}$. Then there is $p\in \mathbb{P}_\alpha$ such that for every $\beta<\gamma$, $p\restriction\gamma+1=p_\beta\restriction\gamma+1$ and $p_\beta\leq p$.
\end{lemma}
\begin{lemma}
Let $\alpha$ be a limit point of $E$ such that $\alpha $ is Mahlo. Then $\mathbb{P}_\alpha$ is $\alpha$-cc.
\end{lemma}

\begin{lemma}
For every $p\in\mathbb{P}_\alpha$ and every statement in the forcing language $\sigma$ there is $p\leq^* p^*$ such that $p^* ||\sigma$.
\end{lemma}
\begin{definition}\label{definition of ultrafilter extension}
 Fix some well ordering $W$ of $V_\lambda$ for some very large $\lambda$ such that for every $\beta<\lambda$ inaccessible $W\restriction V_{\beta}\leftrightarrow \beta$.
 Suppose that $\alpha$ is a limit point of $E$ and $\alpha'=\min\{\beta\in E\mid \beta\geq\alpha\}$ and let $U$ be a normal measure over $\alpha'$. Let us define an ultrafilter $U^*\in V^{\mathbb{P}_\alpha}$:
 \begin{enumerate}
     \item If $\alpha<\alpha'$ then $U^*:=\{X\subseteq\alpha\mid \exists Y\in U, Y\subseteq X\}$.
     \item If $\alpha=\alpha'$ then $E\cap \alpha\notin U$, consider $j_U:V\rightarrow M_U$, since $\alpha\notin j(E)$ it follows that $j_U(\mathbb{P}_\alpha)=\mathbb{P}_\alpha*\mathbb{P}_{(\alpha,j_U(\alpha))}$. Let $G_\alpha$ be $V$-generic. The model $M_U[G_\alpha]$ is closed under $\alpha$-sequences from $V[G_\alpha]$ and let $\l \lusim{A}_\alpha\mid \alpha<\kappa^+\r$ be the $j_U(W)$-minimal enumeration of all the nice names for subsets of $\alpha$. Since the forcing
     $\mathbb{P}_{(\alpha,j_U(\alpha))}/G_\alpha$ has more than $\alpha^+$-closure degree for $\leq^*$ with respect to $V[G_\alpha]$, we can construct a \textit{master sequence}  $\l p_\nu\mid\nu<\alpha^+\r$ such that for each $\nu<\alpha^+$, $p_\nu || \alpha\in j_U(\lusim{A}_\alpha)$. Define $$U^*:=\{(\lusim{A}_\nu)_{G_\alpha}\mid p_\nu\Vdash \alpha\in j_U(\lusim{A}_\alpha)\}$$
 \end{enumerate}
\end{definition}

By \cite{584c9093529345ed9ca7a01a0e0c9d9a}, $U^*$ is a normal ultrafilter over $\alpha'$.
In particular this gives rise to the definition of the normal ultrafilter $U^*_0$ which extends $U_0$ after forcing with $\mathbb{P}_\kappa$. Now we extend $U^*_1$ to a non normal ultrafilter:
Let $j_{U_1}:V\rightarrow M_{U_1}$ be the ultrapower by $U_1$. Then $\kappa\in j_{U_1}(E)$ and therefore $$j_{U_1}(\mathbb{P}_\kappa)=\mathbb{P}_\kappa*Q_\kappa*\mathbb{P}_{(\kappa,j_{U_1}(\kappa))}$$
Where $Q_\kappa$ is the Prikry forcing with the ultrafilter $U'_0$ where $U'_0$ was extended in the same fashion as in Definition \ref{definition of ultrafilter extension} using the minimal witness with respect to $j_0(j_1(W))$. By \cite[Lemma 2.1]{584c9093529345ed9ca7a01a0e0c9d9a}, $U'_0=U^*_0$. Let $\l \lusim{A}_\alpha\mid \alpha<\kappa^+\r$ be the $j_1(W)$-minimal enumeration of all the nice $\mathbb{P}_\kappa$-names for subsets of $\kappa$ and let us define a sequence $\l p_\alpha\mid \alpha<\kappa^+\r$ of $\mathbb{P}_{\kappa+1}$ names for conditions in $\mathbb{P}_{j_{U_1}}(\kappa)/\mathbb{P}_{\kappa+1}$ for a master sequence with respect to $j_1(W)$.
In $V[G]$ define $A\in U^*_1$ if and only if $\exists p\in G\exists B\in U^*_0\exists\nu<\kappa^+$ such that $$p^{\smallfrown}\l\l\r,\lusim{B}\r^{\smallfrown}p_\nu\Vdash\kappa\in j_{U_1}(\lusim{A})$$
where $\lusim{B},\lusim{A}$ are any names interpreted by $G$ to be $B,A$ respectively. By \cite[Lemma 2.2]{584c9093529345ed9ca7a01a0e0c9d9a}, $U^*_1$ is a $\kappa$-complete ultrafilter over $\kappa$ which extends $U_1$. 

The intuition here is that $U^*_1$ concentrates on $\alpha$'s for which we have preformed the Prikry forcing using the $U_0(\alpha)$, hence $[id]_{U^*_1}$ has a Prikry sequence generated. 

Indeed, consider $j_{U^*_1}:V[G]\rightarrow M_{U^*_1}$. Then $j_{U^*_1}(G)$ has a Prikry sequence for each $\alpha\in j_{U^*_1}(E)$. Since $E\in U_1$ and $U_1\subseteq U^*_1$, it follows that $[id]_{U^*_1}$ has a Prikry sequence in $j_{U^*_1}(G)$.
Denote by $\l\kappa_n\mid n<\omega\r$ the Prikry sequence for $[id]_{U^*_1}$. 

Denote $V^*=V[G]$ and consider for each $0<n<\omega$ the function $\psi_n:\kappa\rightarrow [\kappa]^{n}$ defined by $\psi_n(\alpha)=\l\alpha_0,...,\alpha_{n-1}\r$ where $\l\alpha_0,...,\alpha_{n-1}\r$ are the first $n$-elements of the Prikry sequence for $\alpha$ in $G$. The function $\psi_n$ witnesses the Rudin-Keisler projection of $U^*_1$ onto the product of $n$-copies of $U^*_0$ denoted by $U^{*n}_0$:
\begin{proposition}
For each $n<\omega$ and every $X\subseteq[\kappa]^{n}$, $X\in U^{*n}_0$ if and only if $\psi_n^{-1}{}''X\in U^*_1$.
\end{proposition}
\begin{proof}
 Let $X\in U^{*n}_0$, by normality of $U^{*n}_0$, there is $A\in U^*_0$ such that $[A]^n\subseteq X$. Denote by $B:=\psi^{-1''}_n[A]^n=\{\alpha<\kappa\mid \alpha_0,...,\alpha_{n-1}\in A\}$ then there is $p\in G$, $Y\in U^*_0$ and $\nu<\kappa^+$ such that $$p^{\smallfrown}\l\l\r,\lusim{Y}
 \r^{\smallfrown}p_\nu || \kappa\in j_{U_1}(\lusim{B})$$
 Toward a contradiction suppose that 
 $$p^{\smallfrown}\l\l\r,\lusim{Y}\r^{\smallfrown}p_\nu \Vdash \kappa\notin j_{U_1}(\lusim{B})$$
 Namely, 
 $$p^{\smallfrown}\l\l\r,\lusim{Y}\r^{\smallfrown}p_\nu \Vdash \lusim{\kappa}_0,...,\lusim{\kappa}_{n-1}\notin j_{U_1}(\lusim{A})\cap\kappa=\lusim{A}$$
 To see the contradiction, note that $A\in U^*_0$, hence $$p^{\smallfrown}\l\l\r,\lusim{Y}\cap\lusim{A}\r^{\smallfrown} p_\nu \Vdash \lusim{\kappa}_0,...,\lusim{\kappa}_{n-1}\notin\lusim{A}.$$

 As for the other direction, suppose that $X\subseteq[\kappa]^n$ and $B:=\psi^{-1''}_nX\in U^*_1$, then there are $p\in G$, $Y\in U^*_0$ and $\nu<\kappa^+$ such that 
 $$p^{\smallfrown}\l\l\r,\lusim{Y}\r^{\smallfrown} p_\nu \Vdash \kappa\in j_{U_1}(\lusim{B})$$

 It must be the case that for some extension $q\in G$ of $p$, $q\Vdash \lusim{Y}^n\subseteq \lusim{X}$ since otherwise, there would be $\vec{\nu}$ and $p\leq q\in G$ such that $q\Vdash \vec{\nu}\in \lusim{Y}^n\setminus \lusim{X}$. Then $q^{\smallfrown}\l\vec{\nu},\lusim{Y}\setminus\max(\nu)\r^{\smallfrown} p_\nu$ forces that $\vec{\nu}=\l\lusim{\kappa}_0,...,\lusim{\kappa}_{n-1}\r$ and by definition of $B$, $\l\lusim{\kappa}_0,...,\lusim{\kappa}_{n-1}\r\in \lusim{X}$, which is a contradiction. Hence $Y^n\subseteq X$, concluding that $X\in U^{*n}_0$
\end{proof}

Denote by $\l j^*_{n,m}:M^*_n\rightarrow M^*_m\mid n\leq m<\omega\r$ the $\omega$-th iteration by $U^*_0$.  Then by the previous proposition there is a factor map $k_n:M^*_n\rightarrow M_{U^*_1}$ derived by the projection $\psi_n$ i.e. $k_n([f]_{U^{*n}_0})=[f\circ\psi_n]_{U^*_1}$.

\begin{proposition}
For every $n<\omega$, $crit(k_n)=j^*_{0,n}(\kappa)=\kappa_n$, where $\kappa_n$ is the $n$-th element in the Prikry sequence for $[id]_{U^*_1}$.
\end{proposition}
\begin{proof}
 For every $n<\omega$, $\l \kappa, j_{0,1}(\kappa),...,j_{0,n-1}(\kappa)\r=[id]_{U^{*n}_0}$ hence $$k_n(\l \kappa, j_{0,1}(\kappa),...,j_{0,n-1}(\kappa)\r)=[id\circ\psi_n]_{U^{*}_1}=\l\kappa_0,...,\kappa_{n-1}\r.$$ Therefore for every $n<m<\omega$, $j_{0,n}(\kappa)\leq k_m(j_{0,n}(\kappa))=\kappa_n$. We need to prove two separate statements, first that $crit(k_n)=\kappa_n$ and second that $crit(k_n)=j_{0,n}(\kappa)$. Once we prove that $crit(k_m)=\kappa_m$, we can deduce that $j_{0,n}(\kappa)\leq \kappa_n<\kappa_m=crit(k_m)$ which in turn implies that $j_{0,n}(\kappa)=k_m(j_{0,n}(\kappa))=\kappa_n=crit(k_n)$.
 
 Let $\gamma<crit(k_m)$, then $\gamma=[f]_{U^{*n}_0}$ and $\gamma=k_m(\gamma)=[f\circ\psi_m]_{U^*_1}$. 
 Toward a contradiction, suppose that $\gamma\geq \kappa_m$. Thus $$B=\{\alpha<\kappa\mid f(\alpha_0,...,\alpha_{m-1})\geq \alpha_m\}\in U^*_1$$ which by definition implies that there are $p\in G$, $X\in U^*_0$ and $\nu<\kappa^+$ such that $$p^{\smallfrown}\l\l\r,\lusim{X}\r^{\smallfrown} p_\nu\Vdash \kappa \in j_{U_1}(\lusim{B}).$$
 By definition of $B$ it follows that
 $$p^{\smallfrown}\l\l\r,\lusim{X}\r^{\smallfrown} p_\nu\Vdash \lusim{f}(\lusim{\kappa}_0,...,\lusim{\kappa}_{m-1})=j_{U_1}(\lusim{f})(\lusim{\kappa}_0,...,\lusim{\kappa}_{m-1})\geq\lusim{\kappa}_m.$$
In $V[G]$ consider $C_f=\{\alpha<\kappa\mid f''[\alpha]^n\subseteq \alpha\}$. Then $C_f$ is a club and by normality $C_f\in U^*_0$. The condition
 $$p^{\smallfrown}\l\l\r,\lusim{X}\cap\lusim{C}_f\r^{\smallfrown} p_\nu$$
forces that $\lusim{\kappa}_m\in \lusim{C_f}$ which implies that $\lusim{f}(\lusim{\kappa}_0,...,\lusim{\kappa}_{m-1})<\lusim{\kappa}_m$, contradiction.

As for the other direction, let us prove that $\kappa_m\subseteq Im(k_m)$. Let $[f]_{U^*_1}<\kappa_m$, then $B=\{\alpha<\kappa\mid f(\alpha)<\alpha_m\}\in U^*_1$. Hence there are $p\in G$, $X\in U^*_0$ and $\nu<\kappa^+$ such that $$p^{\smallfrown}\l\l\r,\lusim{X}\r^{\smallfrown} p_\nu\Vdash \kappa\in j_{U_1}(\lusim{B}).$$

 Work in $V[G]$, for every $\vec{\gamma}=\l\gamma_0,...,\gamma_{m-1},\gamma_m\r\in [X]^{m+1}$, 
 $$\l\vec{\gamma},X\setminus\max(\vec{\gamma})\r^{\smallfrown} p_\nu\Vdash j_{U_1}(\lusim{f})(\kappa)<\gamma_m.$$
 Hence, by the Prikry property, for each $\delta<\gamma_m$ there is $\nu_\delta<\kappa^+$ and $B_\delta$ such that $\l\vec{\gamma},B_\delta\r^{\smallfrown} p_{\nu_\delta}\Vdash j_{U_1}(\lusim{f})(\kappa)=\delta$. Let
 $$X(\vec{\gamma})=\bigcap_{\delta<\gamma_m}B_\delta\in U^*_0\text{ and } \nu(\vec{\gamma})=\sup_{\delta<\gamma_m}\nu_\delta<\kappa^+.$$ It follows that for some value $\alpha(\vec{\gamma})<\gamma_m$, $$\l \vec{\gamma}, X(\vec{\gamma})\r^{\smallfrown}p_{\nu(\vec{\gamma})}\Vdash j_{U_1}(\lusim{f})(\kappa)=\alpha(\vec{\gamma})<\gamma_m$$
 By varying $\gamma_m$ and pressing down, we can assume that $\alpha(\l\gamma_0,...,\gamma_{m-1},\gamma_m\r)$ does not depend on the last coordinate and hence we can write $$\alpha(\l\gamma_0,...,\gamma_{m-1},\gamma_m\r)=\beta(\l\gamma_0,...,\gamma_{m-1}\r).$$ By normality and regularity respectively, $$X^*:=\Delta_{\vec{\gamma}\in[\kappa]^{m+1}}X(\vec{\gamma})\in U^*_0, \ \nu^*:=\sup_{\vec{\gamma}\in[\kappa]^{m+1}}\nu(\vec{\gamma})<\kappa^+$$
  In particular, there is $p\leq q\in G$ such that $$q^{\smallfrown}\l\l\r,\lusim{X^*}\r^{\smallfrown}p_{\nu^*}\Vdash   j_{U_1}(\lusim{f})(\kappa)=\lusim{\beta}\circ\lusim{\psi_m}(\kappa)$$
 Hence $[f]_{U^*_1}=[\beta\circ\psi_m]_{U^*_1}=k_m([\beta]_{U^{*n}_0})$.

\end{proof}
\section{The failure of Galvin's property for filters}\label{SectionNonGalvinFilters}
Galvin's theorem asserts that if $2^\kappa=\kappa^+$ and $\mathscr{F}$ is a normal filter over $\kappa^+$ then ${\rm Gal}(\mathscr{F},\kappa^+,\kappa^{++})$ holds. It is  natural to ask to what extent this is sensitive to a modification of the above assumptions. As first noticed by Abraham and Shelah \cite{MR830084}, the requirement $2^\kappa=\kappa^+$ is critical. Specifically, it is consistent with \textsf{ZFC} that $2^\kappa>\kappa^+$ and ${\rm Gal}(\mathrm{Cub}_{\kappa^+},\kappa^+,\kappa^{++})$ fails.  

In  \cite[Question~7.8]{TomMotiII} the authors ask whether normality is also a critical requirement. More precisely, it is asked whether there is a  $\kappa^+$-complete filter $\mathscr{F}$ over $\kappa^+$ for which ${\rm Gal}(\mathscr{F},\kappa^+,\kappa^{++})$ fails, yet  $2^\kappa=\kappa^+$ holds.\footnote{Of course, if such a filter exists it cannot be normal.} Our goal in this section is to give a positive answer to this question. The key observation here is that the existence of a $\kappa$-independent family $\mathcal{F}\s\mathcal{P}(\kappa)$ gives rise to such filters, and every \emph{strongly regular cardinal} $\kappa$ (i.e., $\kappa^{<\kappa}=\kappa$) carries such a family. Later we will prove the consistency of every regular cardinal $\kappa$ carrying a $\kappa$-independent family without effecting the power-set-function pattern (see Theorem~\ref{ManyFailuresWithoutTheclubfilter}). 
It is worth to stress that this result does not require any large cardinal  whatsoever. Nevertheless, it is not clear that this construction provides filters containing the club filter. This issue, which refers to \cite[Question~7.9]{TomMotiII},  will be addressed  in \S\ref{SectionManyNonGalvinWithClubfilter}. Specifically, we will show that it is consistent with the \textsf{GCH} that the successor of every singular cardinal $\kappa^+$ carries a non-Galvin filter extending $\mathrm{Cub}_{\kappa^+}$ (see Theorem~\ref{TheoremManyFailures}).

\subsection{Many non-Galvin filters}\label{SectionManyNonGalvin} 
Let $\kappa$ be a regular cardinal. Recall that a family $\mathcal{F}\s\mathcal{P}(\kappa)$ is called \emph{$\kappa$-independent} if for every disjoint subfamilies $\mathcal{F}_1,\mathcal{F}_2\in[\mathcal{F}]^{<\kappa}$ the following is true: $$(\bigcap_{X\in\mathcal{F}_1}X)\cap(\bigcap_{Y\in\mathcal{F}_2}(\kappa\setminus Y))\neq\emptyset.$$\label{independentfamily}
One can show that if $\kappa^{<\kappa}=\kappa$ there is a $\kappa$-independent family of size $2^\kappa$ (see \cite[Exercise~8.10]{Kunen}).
It should be noted that a $\kappa$-independent family does satisfy the $<\kappa$-intersection property, and therefore the set $$\mathscr{W}_{\mathcal{F}}:=\{X\subseteq\kappa\mid \exists\mathcal{A}\in[\mathcal{F}]^{<\kappa}\,(\bigcap\mathcal{A}\subseteq X)\}$$
defines a  $\kappa$-complete filter which includes $\mathcal{F}$. Actually, it is the minimal such filter.
The connection between this concept and non-Galvin filters is exemplified by the next theorem:
\begin{lemma}
    \label{Independent family}
If $\mathcal{F}$ is a $\kappa$-independent family then $\mathscr{W}_\mathcal{F}$ is a $\kappa$-complete filter over $\kappa$ and $\mathcal{F}$ witnesses the failure of  $\mathrm{Gal}(\mathscr{W}_{\mathcal{F}},\kappa,|\mathcal{F}|)$.
\end{lemma}
\begin{proof}
Assume toward contradiction that $\mathcal{F}$ is a $\kappa$-independent family and $\mathscr{W}_{\mathcal{F}}$ is Galvin. Then, there is $\mathcal{A}\in[ \mathcal{F}]^{\kappa}$  such that $\bigcap\mathcal{A}\in \mathscr{W}_{\mathcal{F}}$. By definition of $\mathscr{W}_{\mathcal{F}}$, there is  $\mathcal{B}\in[\mathcal{F}]^{<\kappa}$ such that $\bigcap\mathcal{B}\subseteq \bigcap\mathcal{A}$. Since $|\mathcal{A}|=\kappa$ and $|\mathcal{B}|<\kappa$, one can pick $Y\in\mathcal{A}\setminus\mathcal{B}$. By $\kappa$-independence,  there is $\nu\in \bigcap\mathcal{B}$ such that $\nu\notin Y$ hence, in particular,  $\nu\notin \bigcap\mathcal{A}$. This is a contradiction with the inclusion $\bigcap\mathcal{B}\subseteq \bigcap\mathcal{A}$.
\end{proof}
\begin{theorem}\label{ManyFailuresWithoutTheclubfilter}
   There is a generic extension where: 
   \begin{enumerate}
       \item $\mathsf{GCH}$ holds;
       \item for every regular cardinal $\kappa$ there is a $\kappa$-complete filter $\mathscr{F}$ over $\kappa$ such that $\mathrm{Gal}(\mathscr{F},\kappa,\kappa^+)$ fails.
   \end{enumerate}
\end{theorem}
\begin{proof}
Over $L$ force with $\mathbb{A}$, the class-Easton-supported iteration forcing with $\Add(\kappa,\kappa^+)$ at each regular cardinal $\kappa$. Clearly, this forcing preserves the \textsf{GCH}-pattern of the ground model. For each regular cardinal $\kappa$, $\mathbb{A}_{\kappa+1}$ forces a $\kappa$-independent family of size $\kappa^+$: Let $\langle c_\alpha\mid \alpha<\kappa^+\rangle$ be the Cohen generics introduced by $\Add(\kappa,\kappa^+)$ over $V^{\mathbb{A}_\kappa}$. Next, define $$\mathcal{F}:=\{c^{-1}_\alpha\{1\}\mid \alpha<\kappa^+\}.$$
Using genericity, it is possible to show that $\mathcal{F}$ is $\kappa$-independent. In addition,   the tail forcing $\mathbb{A}/\mathbb{A}_{\kappa+1}$ does preserve this fact, for it defines a $\kappa^+$-directed-closed poset. Thus, the purported result follows.
\end{proof}

	\begin{remark}
The above argument easily adapts to the case where instead of $\Add(\kappa,\kappa^{+})$ we force with $\Add(\kappa,E(\kappa))$, where $E\colon\mathrm{Reg}\rightarrow \mathrm{Card}$ is an Easton class-function. In particular, we can replace the \textsf{GCH} by any allowed power-set-function pattern for regular cardinals.
	\end{remark}

Note that the filters described above do not necessarily extend the club filter, this seems to be quite a restrictive requirement. By the results in \cite{OnPrikryandCohen}, such a filter (ultrafilter) can be forced together with $\textsf{GCH}$ starting from just a single measurable cardinal. If one wishes to get such a filter on a successor cardinal, the cardinals below the measurable cardinal $\kappa$ can be collapsed in order to make $\kappa$ a successor cardinal. It is possible to prove now that the the former ultrafilter which witnessed the failure of the Galvin property before the collapses, generates a non-Galvin filter which extends the club filter. However, it is not clear weather large cardinals are even needed in this situation: 
\begin{question}
\label{qstrgch} What is the consistency strength of $\neg{\rm Gal}(\mathscr{F},\kappa^+,\kappa^{++})$ where $\mathscr{F}$ is a $\kappa^+$-complete filter  over $\kappa^+$ extending $\mathrm{Cub}_{\kappa^+}$ and $2^\kappa=\kappa^+$? Does it require large cardinals?
\end{question}

\subsection{Many non-Galvin filters at successors of singulars containing the club filter}\label{SectionManyNonGalvinWithClubfilter}
The aim of this section is to prove the following result:
\begin{theorem}\label{TheoremManyFailures}
	Assume the $\mathsf{GCH}$ holds and that there is a $\kappa^{+3}$-supercompact cardinal $\kappa$. Then, there is a model of $\mathsf{ZFC}$ where: 
	\begin{enumerate}
		\item $\mathsf{GCH}$ holds; 
		\item for every $\xi\in \mathrm{Ord}$ there is  an $\aleph_{3\cdot\xi+1}$-complete filter $\mathscr{F}_\xi$ over $\aleph_{3\cdot\xi+1}$ that is  is not normal (yet extends the club filter $\mathrm{Cub}_{\aleph_{3\cdot\xi+1}}$) and for which $\mathrm{Gal}(\mathscr{F}_{\xi},\aleph_{3\cdot\xi+1}, \aleph_{3\cdot\xi+2})$ fails. 
	\end{enumerate}
	In particular, the $\mathsf{GCH}$ is consistent with every successor of a singular cardinal $\kappa^+$ carrying a non-normal $\kappa^+$-complete filter $\mathscr{F}$ such that $\mathrm{Cub}_{\kappa^+}\s \mathscr{F}$ and     $\mathrm{Gal}(\mathscr{F},\kappa^+,\kappa^{++})$ fails.
	\end{theorem}
	
	One of the key  tools to produce the above configuration is Abraham and Shelah poset to force the failure of Galvin's property:

\begin{definition}[Abraham-Shelah forcing {\cite{MR830084}}]\label{AbSh forcing}
Let $\kappa<\lambda$ be cardinals with $\kappa$ regular. 
The following clauses yield the definition of the poset $\mathbb{S}(\kappa,\lambda)$:
\begin{itemize}
    
    \item $\lusim{\mathbb{R}}$ is a $\Add(\kappa,1)$-name for the forcing which consists of closed bounded sets $C\subseteq\kappa^+$ which do not contain a subset of cardinality $\kappa$ from $V$.
    The order of $\lusim{\mathbb{R}}$ is forced to be end-extension.
    \item Denote $\mathbb{S}=\Add(\kappa,1)*\lusim{\mathbb{R}}$.
    \item $\mathbb{S}(\kappa,\lambda)$ is a product of $\lambda$-many copies of $\mathbb{S}$ with mixed support, ${<}\kappa$-support on the $\Add(\kappa,1)$-side and $\kappa$-support on the $\lusim{\mathbb{R}}$-side.
    \end{itemize}
\end{definition}
The poset $\mathbb{S}(\kappa,\lambda)$ is $\kappa$-directed closed. Also, assuming that $2^\kappa=\kappa^+$, $\mathbb{S}(\kappa,\lambda)$ is a $\kappa^{++}$-cc forcing that moreover preserves $\kappa^+$ (see \cite[\S1]{MR830084}). This latter fact is a consequence of the following key lemma:
\begin{lemma}[{\cite[\S1]{MR830084}}]\label{Distributivityofthequotient}
    $\mathbb{S}(\kappa,\lambda)$ projects onto $\Add(\kappa,\lambda)$ and the quotient forcing is $(\kappa^+,\infty)$-distributive.
\end{lemma}

Let us now tackle the proof of Theorem~\ref{TheoremManyFailures}. Assume the \textsf{GCH} holds. Let $\mathbb{P}$ be the Easton-supported iteration forcing with $\mathbb{S}(\alpha^+,\alpha^{+3})\ast \lusim{\col}(\alpha,\alpha^{+})$ at every inaccessible cardinal $\alpha\leq \kappa.$ It is easy to show that forcing with $\mathbb{P}$ preserves the \textsf{GCH} pattern. Also, standard arguments show that this poset preserves the $\kappa^{+3}$-supercompactness of $\kappa$.\footnote{The reason for this being that the target model of a $\kappa^{+3}$-supercompact embedding is closed enough to cook up a master condition.}
Let $V$ denote the generic extension by $\mathbb{P}_\kappa\ast\lusim{\mathbb{S}}(\kappa^+,\kappa^{+3})$. Observe that in this model $2^\kappa=\kappa^+$ holds while  $\Gal{\kappa^{++}}{\kappa^{+3}}$ fails.
			
			\begin{lemma}\label{TheFilter}
			Working in $V^{\col(\kappa,\kappa^+)}$, there is a $\kappa^+$-complete filter $\mathscr{F}$ over $\kappa^+$ that extends $\mathrm{Cub}_{\kappa^+}$ and witnesses that  $\mathrm{Gal}(\mathscr{F},\kappa^+,\kappa^{++})$ fails.
			\end{lemma} 
			\begin{proof}
			Working in $V^{Col(\kappa,\kappa^+)}$ note that $\kappa^+=(\kappa^{++})^V$ and let $\mathscr{F}$ be the filter generated by $(\mathrm{Cub}_{\kappa^{++}})^V$. Namely,
			$$\mathscr{F}:=\{X\s\kappa^{+}\mid \exists C\in\mathrm{Cub}_{\kappa^{++}}^V\, (C\s X)\}.$$
			Since $2^\kappa=\kappa^+$, $\col(\kappa,\kappa^+)$ is $\kappa^{++}$-cc. In particular, every $C\in\mathrm{Cub}_{\kappa^+}^{V^{\col(\kappa,\kappa^+)}}$ contains some $D\in\mathrm{Cub}_{\kappa^{++}}^V$. 
			It thus follows that $\mathrm{Cub}_{\kappa^+}^{V^{\col(\kappa,\kappa^+)}}\s \mathscr{F}$. 
			
			\smallskip
			

			We next show that $\mathscr{F}$ is $\kappa^+$-complete and  $\mathrm{Gal}(\mathscr{F},\kappa^+,\kappa^{++})$ fails.
			
			\begin{claim}\label{PreservingFailure}
			$\mathrm{Gal}(\mathscr{F},\kappa^+,\kappa^{++})$ fails. 	
			\end{claim}
			\begin{proof}[Proof of claim]
				Let $\langle X_\alpha\mid \alpha<\kappa^{+3}_V\rangle\in V$  be a sequence witnessing the failure of $\mathrm{Gal}(\mathrm{Cub}_{\kappa^{++}},\kappa^{++},\kappa^{+3})$. Let  $C$ be a ${\col(\kappa,\kappa^+)}$-generic over $V$. We claim that the above sequence witnesses $\neg \mathrm{Gal}(\mathscr{F},\kappa^+,\kappa^{++})$ in $V[C]$. Otherwise, let $I\in[\kappa^{++}]^{\kappa^+}\cap V[C]$  such that $\bigcap_{\alpha\in I} X_\alpha\in \mathscr{F}$. By definition, there is a set $Y_I\in\mathrm{Cub}^V_{\kappa^{++}}$ contained in this intersection. By our assumption, $\{\alpha<\kappa^{+3}_V\mid Y_I\s X_\alpha\}$ has cardinality ${\leq}\kappa^{+}$ in $V$, hence cardinality ${\leq}\kappa$ in  $V[C]$. However, $I\s\{\alpha<\kappa^{++}\mid Y_I\s X_\alpha\}$, which yields a contradiction.
			\end{proof}

			\begin{claim}\label{Preservingcompleteness}
			Let $\mathbb{P}$	 be a $\theta$-cc forcing notion and $\mathscr{W}$ a $\theta$-complete filter over a set $I$. Then, $\mathscr{F}_{\mathscr{W}}:=\{X\s I\mid \exists W\in\mathscr{W}\, (W\s X)\}$ is $\theta$-complete in $V^{\mathbb{P}}$.
			\end{claim}
			\begin{proof}[Proof of claim]
			Let $p\in\mathbb{P}$, $\lambda<\theta$ and $\langle \lusim{X}_\alpha\mid \alpha<\lambda\rangle$ be a sequence of $\mathbb{P}$-names such that $p\forces_{\mathbb{P}}\langle \lusim{X}_\alpha\mid \alpha<\lambda\rangle\s \lusim{\mathscr{F}}_\mathscr{W}.$ For each $\alpha<\lambda$ let $\mathcal{A}_\alpha\s \mathbb{P}/p$ be a maximal antichain such that for each $q\in\mathcal{A}_\alpha$ there is $W_{q,\alpha}\in\mathscr{W}$ with  $q\forces_{\mathbb{P}}\lusim{X}_\alpha\supseteq W_{q,\alpha}$. In particular, $p\forces_{\mathbb{P}} \lusim{X}_\alpha\supseteq W_\alpha$, where $W_\alpha:=\bigcap_{q\in \mathcal{A}_\alpha}W_{q,\alpha}.$
			
			Note, however, that $\langle W_\alpha\mid \alpha<\lambda\rangle$ might not belong to $V$. To work around this we use yet again the $\theta$-ccness of $\mathbb{P}$: Let $\lusim{f}\colon \lambda\rightarrow \mathscr{W}$ be a $\mathbb{P}$-name for the above sequence and find $F\in V$, $F\colon \lambda\rightarrow\mathcal{P}_{<\theta}(\mathscr{W})$ such that $p\forces_{\mathbb{P}}``\forall\alpha<\lambda\,(\lusim{f}(\alpha)\in \check{F}(\alpha))$''. Set, $W:=\bigcap_{\alpha<\lambda}\bigcap F(\alpha)$. Note that $W\in \mathscr{W}$, by $\theta$-completeness of $\mathscr{W}$. Clearly, 
			$p\forces_{\mathbb{P}} W\s \bigcap_{\alpha<\lambda}\lusim{X}_\alpha$. This shows that $p\forces_{\mathbb{P}}\bigcap_{\alpha<\lambda}\lusim{X}_\alpha\in\lusim{\mathscr{F}}_{\mathscr{W}}$, as wanted.
			\end{proof}
			In particular, $\mathscr{F}$ is a $\kappa^+$-complete filter in  $V^{\col(\kappa,\kappa^+)}$. 
\end{proof}

\begin{remark}
We indicate that in the above claim we obtain only the \emph{weak} failure of Galvin's property. Indeed,  we do not know whether the strong failure is forceable in this context.
\end{remark}
As of now we have produced a model of \textsf{GCH} where $\kappa$ is $\kappa^{+3}$-supercompact and there is a $\kappa^+$-complete filter $\mathscr{F}$ over $\kappa^+$ such that $\mathrm{Cub}_{\kappa^+}\s \mathscr{F}$ and $\neg \mathrm{Gal}(\mathscr{F},\kappa^+,\kappa^{++})$. By Galvin's theorem, $\mathscr{F}$ cannot be normal. In a slight abuse of notation we yet again denote this model by $V$. 

\smallskip

Let us agree that $(j,F)$ is a \emph{weak constructing pair} if it satisfies the clauses of \cite[Definition~1.7]{NegGalSing} with the only exception that $\col(\kappa^{+4},i(\kappa))^N$ is replaced by $\col(\kappa^{+3},i(\kappa))^N$. As the \textsf{GCH} holds   it is quite easy to produce such a pair with $j\colon V\rightarrow M$ witnessing $\kappa^{++}$-supercompactness of $\kappa$ (see e.g., \cite[Lemma~2.4]{NegGalSing} or \cite[Lemma~8.5]{Cummings-handbook}). Let $u_*$ be the  measure sequence inferred from $(j,F)$. Arguing  as in \cite[Lemma~1]{CumGCH}, $u_*\restriction\alpha$ exists and belongs to $\mathcal{U}_\infty$, for all 
$\alpha<\kappa^{+++}$. In particular, there is $\alpha<\kappa^{+3}$ such that the sequence $u:=u_*\restriction\alpha\in \mathcal{U}_\infty$ has a repeat point.\footnote{I.e., an ordinal $\gamma<\len(u)$ such that $\mathscr{F}_u=\mathscr{F}_{u\restriction\gamma}$.} Now, let $\mathbb{R}_u$ be the Radin forcing with interleaved collapses as defined in \cite[\S1]{NegGalSing} with the minor change described above. The next lemma can be proved exactly as in \cite[Proposition~1.34]{NegGalSing}:

\begin{lemma}[Cardinal structure]\label{CardinalStructure}
 The following holds in $V^{\mathbb{R}_u}$:
\begin{enumerate}
    \item Every $V$-cardinal $\geq\kappa^{+}$ remians a cardinal;
    \item The only cardinals $\leq \kappa$  are
    $$\{\aleph_0,\aleph_1, \aleph_2, \aleph_3\}\cup \{(\kappa_{\xi}^{+k})^V\mid 1\leq k\leq 3,\, \xi<\kappa\}\cup\Lim(C)\cup\{\kappa_u\};$$
\end{enumerate}
Also, if $u$ has a repeat point then  $\kappa$ remains measurable in $V^{\mathbb{R}_u}$.
\end{lemma}

Since $j\colon V\rightarrow M$ is a   $\kappa^{++}$-supercompact embedding, $\mathrm{Gal}(\mathscr{F},\kappa^+,\kappa^{++})$ fails in the model $M$. Therefore, 
$$A:=\{w\in \mathcal{U}_\infty\mid \text{$``\exists\mathscr{F}_{w}$ witnessing Lemma~\ref{TheFilter} w.r.t. $\kappa^+_w$''}\}\in\mathscr{F}_{u}.$$

Let $G\s\mathbb{R}_u/p$ be $V$-generic, where $p:=\langle (u, \omega, \emptyset, A, H)\rangle$ and $H\in F^*_u$. Everything is now in place to complete the proof of Theorem~\ref{TheoremManyFailures}:
\begin{lemma}
 The following properties are true in $V[G]_\kappa$:
 \begin{enumerate}
 	\item $\mathsf{GCH}$ holds;
 	\item for every $\xi\in\mathrm{Ord}$ there is an $\aleph_{3\cdot\xi+1}$-complete filter $\mathscr{F}_\xi$ over $\aleph_{3\cdot\xi+1}$ that is  is not normal (yet extends the club filter $\mathrm{Cub}_{\aleph_{3\cdot\xi+1}}$) and for which $\mathrm{Gal}(\mathscr{F}_{\xi},\aleph_{3\cdot\xi+1}, \aleph_{3\cdot\xi+2})$ fails. 
  \end{enumerate}
\end{lemma}
\begin{proof}
	We divide the proof into a series of claims:
	\begin{claim}
		Clause~(1) holds.
	\end{claim}
	\begin{proof}[Proof of claim]
	This follows from an easy counting-nice-name argument involving the \textsf{GCH} from $V$ and the usual factoring of Radin forcing. 
	\end{proof}
	
	\begin{claim}
		Clause~(2) holds.
	\end{claim}
	\begin{proof}[Proof of claim]
Let $\langle \kappa_\xi\mid \xi<\kappa\rangle$ be the increasing enumeration of the generic club induced by $G$.  By Proposition~\ref{CardinalStructure}, $\kappa$ remains inaccessible in $V[G]$. 

For each  $\xi<\kappa$ write
$$\Phi(\xi)\equiv \text{``$\mathscr{F}_{\xi}$ witnesses Lemma~\ref{TheFilter}''},$$
where $\mathscr{F}_\xi$ is a filter witnessing $u_\xi\in A$.  

Fix $\xi<\kappa$ and pick $q\in G$ mentioning both $\kappa_\xi$ and $\kappa_{\xi+1}$, say at coordinates $m$ and $m+1$, respectively. 
Note that $\Phi(\xi)$ holds in $V$ by virtue of our choice of $A$. The usual factoring arguments give:
$$\mathbb{R}_u/q\simeq \mathbb{R}_{u_{\xi}}/q^{\leq m}\times \col(\kappa_{\xi}^{+3},\kappa_{{\xi+1}})\times\mathbb{R}_u/q^{>m+1}.$$
To not complicate the notations we shall tend to identify $\mathscr{F}_\xi$ with the filter generated by it in the different sub-generic extensions of $V^{\mathbb{R}_u}$.

\smallskip

$\br$ \underline{$\mathscr{F}_{\xi}$ is $\kappa^+_\xi$-complete:} The first forcing is $\kappa^+_\xi$-cc, hence by Claim~\ref{Preservingcompleteness} it preserves $\kappa^+_\xi$-completeness of $\mathscr{F}_\xi$. The rest does not introduce $\kappa^+_\xi$-sequences so that it also preserves the property under consideration. 

\smallskip

$\br$ \underline{$\mathscr{F}_{\xi}$ extends the club filter:}   Let $C\in\mathrm{Cub}_{\kappa^+_{\xi}}$ in $V[G]$. By the above factoring $C\in V[G_0]$, where $G_0$ is the projection of $G$ on $\mathbb{R}_{u_\xi}/q^{\leq m}$. Since this forcing is $\kappa^{+}_\xi$-cc there is $D\in \mathrm{Cub}_{\kappa^+_\xi}^V$ such that $D\s C$. Finally, since $\mathrm{Cub}_{\kappa^+_\xi}^V\s \mathscr{F}_\xi$ it follows that the filter generated by $\mathscr{F}_\xi$ has $C$ as an element.

\smallskip

$\br$ \underline{$\neg \mathrm{Gal}(\mathscr{F}_\xi,\kappa^+_\xi,\kappa^{++}_\xi)$:} Let $\mathcal{X}=\langle X_\alpha\mid \alpha<\kappa^{++}_\xi\rangle\s \mathscr{F}_\xi$ in $V$  witnessing the failure of $\mathrm{Gal}(\mathscr{F}_\xi,\kappa^+_\xi,\kappa^{++}_\xi)$. 
Since $\kappa_\xi^+$ and $\kappa^{++}_\xi$ are preserved in $V[G_0]$ the argument of Claim~\ref{PreservingFailure} shows that $\mathcal{X}$ is still a witness for the failure of Galvin's property in $V[G_0]$. The other two posets do not introduce new subsets to $\kappa^{++}_\xi$ so $\mathcal{X}$ is going to be a witness for $\neg \mathrm{Gal}(\mathscr{F}_\xi,\kappa^+_\xi,\kappa^{++}_\xi)$ in $V[G]$.



\smallskip

Altogether, $\Phi(\xi)$ holds in $V[G]$. 
	\end{proof}
	Looking at Proposition~\ref{CardinalStructure}(2) it is easy to check that 
$$\kappa^+_\xi:=\begin{cases}
\aleph_{3\cdot(\xi+1)+1}, & \text{if $\xi<\omega$;}\\
\aleph_{3\cdot \xi+1}, & \text{if $\xi\geq \omega$.}
\end{cases}$$
By further forcing with $\col(\aleph_0,\aleph_3)$ we get that, for every ordinal $\xi<\kappa$,  $\kappa_\xi^+=\aleph_{3\cdot\xi+1}$ holds in the resulting generic extension. Moreover, this forcing is $\aleph_4$-cc and, as a result, preserves $\Phi(\xi)$ for every $\xi<\kappa$ (see Lemma~\ref{TheFilter}).
\smallskip
Finally,  $V[G]_\kappa$ has the desired properties. 
\end{proof}
In the light of Theorem~\ref{ManyFailuresWithoutTheclubfilter} and \ref{TheoremManyFailures} the following becomes natural:
\begin{question}
	What is the consistency strength of the configuration described in Theorem~\ref{TheoremManyFailures}? Does it require large cardinals?
	\end{question}

 \section{Failure of Galvin's property for the club filter}\label{SectionFailurefortheClubfilter}
 In this section we continue the study of \cite{MR830084, NegGalSing} and analyze the failure of Galvin's property relative to the club filter. The original motivation for this section was to produce infinitely-many consecutive failures of Galvin's property; e.g., $\neg\Gal{\aleph_{n+1}}{\aleph_{n+2}}$ for all $n<\omega$. Regrettably, we did not succeed in this enterprise (see Question~\ref{QuestionAboutManyfailures}). Instead, we will be analyzing the following failure of \emph{pseudo-compactness}: $\Gal{\kappa^+}{\kappa^{++}}$ holds, yet $\Gal{\alpha^+}{\alpha^{++}}$ fails  for \emph{many} cardinals $\alpha<\kappa$.
 
 In \S\ref{MeasureConcentratingOnfailures} we describe how to force a (non normal) measure $\mathscr{U}$  such that for  $\mathscr{U}$-many $\alpha$'s, $\Gal{\alpha^+}{\alpha^{++}}$ fails. In particular, after forcing with the \emph{Tree Prikry forcing} relative to this measure, one can produce a model where $\Gal{\alpha^+}{\alpha^{++}}$ fails for unboundedly many cardinals below a strong limit cardinal of countable cofinality. Later, in \S\ref{SectionTheModifiedPrikry} we obtain a similar result at the level of the very first singular cardinal, $\aleph_\omega$. Both constructions are performed out of sharp assumptions; i.e., from a measurable cardinal.
 \subsection{A measure concentrating on failures of Galvin's property}\label{MeasureConcentratingOnfailures}
Assume that the \textsf{GCH} holds and let $\mathscr{U}$ be a normal measure over $\kappa$. Let $\mathbb{P}_\kappa$ denote the Easton-supported iteration $\l \mathbb{P}_\alpha,\lusim{\mathbb{Q}}_\beta\mid \alpha\leq\kappa,\beta<\kappa\r$ where,  for each inaccessible 
cardinal $\beta<\kappa$, we force with the lottery sum of $\mathbb{S}(\beta,\beta^{++})$ and the trivial poset. As usual, for non-inaccessible $\beta$'s the iteration forces with the trivial poset. Let $G\subseteq \mathbb{P}_\kappa$ be a $V$-generic filter.

Opting for the trivial forcing at the $\kappa^{\mathrm{th}}$-stage  of $j_{\mathscr{U}}(\mathbb{P}_\kappa)$ we can extend $j_{\mathscr{U}}$ to a $V[G]$-definable embedding  $j_1^*:V[G]\rightarrow M[H]$. For this we use our \textsf{GCH}-assumption and the usual lifting arguments.
Next, let $\alpha\in (\kappa,j_{\mathscr{U}}(\kappa))$ be an ordinal such that $H$ opts for $\mathbb{S}(\alpha,\alpha^{++})$. Note that there is some of such $\alpha$'s by density. 
Define
.$$\mathscr{W}:=\{X\subseteq \alpha\mid \alpha\in j_1^*(X)\}$$
Note that $\mathscr{W}$ is a $\kappa$-complete ultrafilter over $\kappa$ concentrating on the collection of all $\beta$'s for which $\Gal{\beta^+}{\beta^{++}}$ fails. In particular, 
\begin{theorem}
Starting from the assumptions of this section, the following is consistent:
$\kappa$ is a strong limit cardinal, $\cf(\kappa)=\omega$,  $\Gal{\kappa^+}{\kappa^{++}}$ holds, but $\Gal{\alpha^+}{\alpha^{++}}$ strongly fails for cofinaly many $\alpha<\kappa$. 
\end{theorem}
\begin{proof}
Let $\mathscr{W}$ be the $\kappa$-complete measure defined above. Force with $\mathbb{P}_T(\mathscr{W})$, the corresponding Tree Prikry forcing. In the resulting extension $2^\kappa=\kappa^+$ holds, hence $\Gal{\kappa^+}{\kappa^{++}}$ also does. In addition, $\Gal{\alpha^+}{\alpha^{++}}$ strongly fails for all the members $\alpha$ in the Prikry sequence: this is because $\mathbb{P}_T(\mathscr{W})$ does not add bounded subsets to $\kappa$ and $\Gal{\alpha^+}{\alpha^{++}}$ strongly fails in $V$.
\end{proof}
By choosing an arbitrary $\alpha$ we lose control of the measure $\mathscr{W}$. For instance, it is unclear whether $\mathscr{W}$ contains the club filter $\mathrm{Cub}_\kappa.$ To work around this we shall look at the second ultrapower by $\mathscr{U}$ and extract from there a different measure. So, put $\mathscr{U}_2:=j_\mathscr{U}(\mathscr{U})$ and look at $j_{\mathscr{U}_2}:M_\mathscr{U}\rightarrow M_{\mathscr{U}_2}.$

Let us extend the embedding $j_{\mathscr{U}_2}$. For this, we choose to force with $\mathbb{S}(j_{\mathscr{U}}(\kappa),j_{\mathscr{U}}(\kappa)^{++})_{M_{\mathscr{U}_2}}$ at the $j_{\mathscr{U}}(\kappa)^{\mathrm{th}}$-stage of the iteration $j_{\mathscr{U}_2}(j_{\mathscr{U}}(\mathbb{P}_\kappa))$. 
By our \textsf{GCH} assumption,  $|j_{\mathscr{U}_2}(j_{\mathscr{U}}(\mathbb{P}_\kappa))|=\kappa^+$. Using this we can construct in $V[G]$ a generic for $j_{\mathscr{U}_2}(j_{\mathscr{U}}(\mathbb{P}_\kappa))$, where at the $j_{\mathscr{U}}(\kappa)^{\mathrm{th}}$-stage the iteration opts for $\mathbb{S}(j_{\mathscr{U}}(\kappa),j_{\mathscr{U}}(\kappa)^{++})_{M_{\mathscr{U}_2}}$. As in the previous argument, we extend  the embedding  $j_{\mathscr{U}_2}\circ j_{\mathscr{U}}$ to another embedding $j^*_2$.
Define,  $$\mathscr{W}:=\{X\s \kappa \mid \kappa_1\in j^*_2(X)\}.$$
Once again, $\mathscr{W}$ defines  a $\kappa$-complete filter over $\kappa$ concentrating on those $\beta$'s for which $\Gal{\beta^{+}}{\beta^{++}}$ fails. In addition, $\mathscr{W}$ contains the club filter $\mathrm{Cub}_\kappa$, hence there are stationarily-many $\beta$'s for which Galvin's property fails. Note, however, that $\mathscr{W}$ is not normal; should this be the case it will entail the failure of $\Gal{\kappa^{+}}{\kappa^{++}}$ and as a result the failure of the \textsf{SCH}.\footnote{Recall that we just assume the existence of a measurable cardinal.} The issue of getting a normal measure $\mathscr{W}$ exhibiting the above non-Galvin-like pattern will be revisited at the beginning of \S\ref{SectionConsistencyStrength}.

 \subsection{A Prikry-type poset forcing the failure of Galvin's property}\label{SectionTheModifiedPrikry}
Assume the \textsf{GCH} holds and let $\kappa$ be a measurable cardinal. Let $\mathcal{U}$ be a normal measure on $\kappa$ and $j\colon V\rightarrow M$ the corresponding ultrapower.
\begin{lemma}
There is $K\in V$ that is $M$-generic for the poset $$(\mathbb{S}(\kappa^{+},\kappa^{+3})\times\col(\kappa^{+3},{<}j(\kappa))^M.$$
\end{lemma}
\begin{proof}
Denote by $\mathbb{Q}$ the above-displayed poset. Note that $\mathbb{Q}$ is $j(\kappa)$-cc in $M$
and, also, every condition in $\mathbb{Q}$ can be identified with a member of $(V_{j(\kappa)})^M$. Since $j(\kappa)$ is inaccessible in $M$, combining these two facts one infers that every maximal antichain (in $M$) for $\mathbb{Q}$ can be regarded as a member of $(V_{j(\kappa)})^M$. In particular, there are at most $|j(\kappa)|^V$-many such objects in $V$. Using our assumption  that $2^\kappa=\kappa^+$ we conclude that there are at most $\kappa^+$-many of such. Now, since $\mathbb{Q}$ is $\kappa^+$-directed-closed in $M$ and $M^\kappa\s M$, so it is in $V$. From altogether we can easily produce the desired $M$-generic filter $K$ by diagonalizing over all the maximal antichains lying in $M$.
\end{proof}
\begin{remark}
Note that the above argument is not available if we replace $\mathbb{S}(\kappa^{+},\kappa^{+3})$ by, e.g., $\mathbb{S}(\kappa,\kappa^{++})$.
\end{remark}
Hereafter we shall denote $\mathbb{Q}:=(\mathbb{S}(\kappa^{+},\kappa^{+3})\times\col(\kappa^{+3},{<}j(\kappa)))^M.$ For a measure one set of inaccessible $A\in \mathcal{U}$ we have that $\mathbb{Q}:=[\rho\in A\mapsto \mathbb{Q}(\rho)]$, where $\mathbb{Q}(\rho)$ stands for $\mathbb{S}(\rho^+,\rho^{+3})\times \col(\rho^{+3},{<}\kappa).$ Likewise, for an inaccessible cardinal $\rho<\rho'\leq \kappa$ we shall denote $\mathbb{Q}(\rho,\rho'):=\mathbb{S}(\rho^+,\rho^{+3})\times\col(\rho^{+3},{<}\rho').$ 
\begin{definition}\label{ModifiedPrikry}
Let $\mathbb{P}$ the set of all conditions of the form $$p=\langle \rho_0,a_0, \dots, \rho_{\ell(p)-1},a_{\ell(p)-1}, A, H\rangle$$
such that the following requirements are met:
\begin{enumerate}
    \item $\rho_0<\cdots<\rho_{n-1}$ are inaccessible cardinals below $\kappa$;
    \item for each $i<\ell(p)$, $a_i\in \mathbb{Q}(\rho_i,\rho_{i+1})$, where we stipulate $\rho_{\ell(p)}:=\kappa.$
    \item $A\in \mathcal{U}$ and $A\s \{\rho<\kappa\mid\text{$\rho$ is inaccessible and $\rho>\supp(a_{\ell(p)-1})$}\}$;\footnote{Here $\supp(a_{i})$ is the smallest $\beta>\rho_{i}$ such that  condition $a_{i}\in V_\beta$.} 
    \item $H$ is a function with $\dom(H)=A$, $H(\rho)\in \mathbb{Q}(\rho)$ and $[H]_\mathcal{U}\in K.$
\end{enumerate}
We shall refer to $\stem(p):=\langle \rho_0, a_0,\cdots, \rho_{\ell(p)-1},a_{\ell(p)-1}\rangle$ and $H$ as the \emph{stem} and \emph{supplier} of $p$, respectively. In addition, given two conditions 
\begin{eqnarray*}
    p=\langle \rho_0,a_0, \dots, \rho_{\ell(p)-1},a_{\ell(p)-1}, A, H\rangle\\
    q=\langle \rho'_0,b_0, \dots, \rho'_{\ell(q)-1},b_{{\ell(q)}-1}, B, L\rangle
\end{eqnarray*}
we shall write $p\geq^* q$ in case $\ell(p)=\ell(q)$ and the following  hold:
\begin{enumerate}
    \item for each $i<\ell(p)$, $\rho_i=\rho'_i$ and $a_i\geq_{\mathbb{Q}(\rho_i,\rho_{i+1})}b_i$;
    \item $A\s B$ and $H(\rho)\geq L(\rho)$ for all $\rho\in A.$
\end{enumerate}
Finally, we shall write $p\geq^{**}q$ iff $p\geq^* q$ and $a_i=b_i$ for all $i<\ell(p).$
\end{definition}

\begin{definition}\label{onepoint}
Given $p\in\mathbb{P}$ as above and $\rho\in A$ define
$$p\cat \rho:=\langle\rho_0,a_0,\cdots, \rho_{\ell(p)-1},a_{\ell(p)-1},\rho, H(\rho), A_\rho, H_\rho\rangle,$$
where $A_\rho:=\{\varrho\in A\mid \varrho>\supp(H(\rho))\}$ and $H_\rho:=H\restriction A_\rho.$

In general, given $\vec\rho\in [A]^{<\omega}$ the sequence $p\cat{ \vec\rho}$ is defined recursively.\footnote{Here, by convention, $p\cat \emptyset:=p.$} 
\end{definition}
\begin{remark}
For every $p\in \mathbb{P}$ and $\rho\in A$, $p\cat \rho$ is a legitimate member of $\mathbb{P}$. In effect,  since $H$ and $H_\varrho$ differ on a $\mathcal{U}$-negligible set, $[H_\rho]_\mathcal{U}=[H]_\mathcal{U}\in K.$ Clearly, the same applies to $p\cat{\vec\rho}$ for $\vec\rho\in [A]^{<\omega}.$
\end{remark}
\begin{definition}
For $p,q\in\mathbb{P}$ write $p\geq q$ iff there is $\vec\rho\in [A^q]^{<\omega}$ such that $p{}^\curvearrowright{\vec\rho}\geq^* q.$ In addition, for a condition $q\in\mathbb{P}$, we write $t\preceq \stem(q)$ if there is $p\geq q$ such that $t=\stem(p).$
\end{definition}
The following is almost immediate.
\begin{lemma}
$\mathbb{P}$ is $\kappa^+$-cc and even $\kappa^+$-Knaster. In particular, in $V^\mathbb{P}\models 2^\kappa=\kappa^+.$
\end{lemma}
\begin{proof}
Let $A\s \mathbb{P}$ be an antichain of size $\kappa^+$. Without loss of generality we may assume that both the length and the stem of the conditions in $A$ are fixed. Note that this is possible in that any stem is a member of $V_\kappa$. Finally, observe that all members of $A$ are compatible: this is thanks to the requirement that the suppliers come from $K$, which is a filter.
\end{proof}

\begin{lemma}\label{factoring}
Let $p=\langle \rho_0,a_0, \dots, \rho_{n-1},a_{n-1}, A, H\rangle\in \mathbb{P}$. Then, there is an isomorphism between $\mathbb{P}/p$ and $\prod_{i<n-1}\mathbb{Q}(\rho_i,\rho_{i+1})\times \mathbb{P}/\langle \rho_{n-1},a_{n-1}, A,H\rangle.$
\end{lemma}
\begin{proof}
Let $\Phi\colon \mathbb{P}/p\rightarrow \prod_{i<n-1}\mathbb{Q}(\rho_i,\rho_{i+1})\times \mathbb{P}/\langle \rho_{n-1},a_{n-1}, A,H\rangle$ be the map
$$
\Phi\colon q\mapsto\langle \langle a^q_0,\ldots, a^q_{n-2}\rangle, \langle \rho_{n-1},a_{n-1},\rho_n^q,a^q_n,\ldots, \rho_{\ell(q)-1},a^q_{\ell(q)-1} A^q, H^q\rangle \rangle.$$
Clearly, $\Phi$ is bijective, and both it and its inverse are order-preserving.
\end{proof}
The upcoming lemma provides the main technical tool. Prior to proving it let us first introduce some useful notation:
\begin{notation}
For $p\in\mathbb{P}$ and $t\preceq \stem(p)$ we denote by  $p+t$ the sequence 
$p+t:=t^\smallfrown \langle A^p_{\max(t)}, H^p_{\max(t)}\rangle.$ It is easy to check that $p+t\in \mathbb{P}.$
\end{notation}
\begin{lemma}
$\mathbb{P}$ has the Prikry property. 
\end{lemma}
\begin{proof}
 Let $p\in\mathbb{P}$ and $\varphi$ be a sentence in the forcing language of $\mathbb{P}$. We shall separate the proof of the lemma into three claims:
\begin{claim}
There is $p\geq^{**}q$ with the following property: for every $q\geq r$ such that $r\parallel \varphi$ then $q+\stem(r)\parallel \varphi.$
\end{claim}
\begin{proof}[Proof of claim]
Let $t\preceq\stem(p)$ be an arbitrary stem. If there is some $r\geq p$ with $\stem(r)=t$ and $r\parallel \varphi$ then put $H_t:=H^r$. Otherwise, put $H_t:=H$. This  generates a directed set of conditions $\{[H_t]_\mathcal{U}\mid t\preceq \stem(p)\}$ in $\mathbb{Q}$, which belongs to $M$.\footnote{Note here the need for closure under $\kappa$-sequences on $M$.}  In particular, we can let $[H^*]_\mathcal{U}$ be a $\leq_{\mathbb{Q}}$-lower bound for this collection. Now, for each $t$ consider
$$A_t:=\{\rho<\kappa\mid H^*(\rho)\geq_{\mathbb{Q}(\rho)} H_t(\rho)\}.$$
Clearly, $A^*:=\dom(H^*)\cap \diagonal_{t\preceq\stem(p)} A_t\in \mathcal{U}.$ Put $q=\stem(p)^\smallfrown \langle A^*,H^*\rangle.$

We claim that $q$ has the desired property. To show this let $r\geq q$ be a condition that decides $\varphi$. Since $r\geq p$ there is a (possibly different) condition $u\geq p$ deciding $\varphi$ with $\stem(u)=\stem(r)$ and such that $H_{\stem(u)}=H^{u}.$ Note that the very definition of diagonal intersection yields $q+\stem(r)\geq u$, hence $q+\stem(r)\parallel\varphi,$ as well. This completes the verification of the claim.
\end{proof}
Let $q=s^\smallfrown \langle A^*, H^*\rangle$ be the condition obtained in the previous claim.
\begin{claim}
There is a condition $q^*\geq^{**} q$ with the following property: 

For each $t\preceq s$, if $\rho\in A^{q^*}_{\max(t)}$ and $a\geq_{\mathbb{Q}(\rho)} H^{q^*}_{\max(t)}(\rho)$ are such that
$t^\smallfrown \langle \rho, a\rangle^\smallfrown \langle A^{q^*}_\rho, H^{q^*}_{\rho}\rangle\parallel \varphi$
then $(q^*+t)\cat \rho\parallel \varphi$ for all $\rho\in A^{q^*}_{\max(t)}.$ 

Moreover, the decision made by the conditions $(q^*+t)\cat \rho$ is uniform.
\end{claim}
\begin{proof}
For a stem $t\preceq s$  consider
$$ 
D^0_t:=\{[L]_\mathcal{U}\in \mathbb{Q}/[H^*]_\mathcal{U}\mid (t^\smallfrown \langle\kappa, [L]_\mathcal{U}\rangle^\smallfrown \langle j(A^*)_\kappa, j(H^*)_\kappa\rangle\parallel j(\varphi))\}
$$
and 
$$ 
D^1_t:=\{[L]_\mathcal{U}\in \mathbb{Q}/[H^*]_\mathcal{U}\mid \forall [L']_\mathcal{U}\geq_{\mathbb{Q}} [L]_\mathcal{U}\;
(t^\smallfrown \langle\kappa, [L']_\mathcal{U}\rangle^\smallfrown \langle j(A^*)_\kappa, j(H^*)_\kappa\rangle\nparallel j(\varphi))\}.
$$
Clearly, $D_t:=D^0_t\cup D^1_t$  is dense below $[H^*]_\mathcal{U}\in K$. Hence, we can choose $[H_t]_\mathcal{U}\in K\cap D_t$. Yet again, this yields a directed set $\{[H_t]_\mathcal{U}\mid t\preceq_\mathbb{Q} s\}$  (in $M$) of conditions in $K$. Let $[\bar{H}]_\mathcal{U}$ be  a $\leq_\mathbb{Q}$-lower bound  for this collection.

For $t\preceq s$ there is $i(t)\in \{0,1\}$ such that $[H_t]_\mathcal{U}\in D^{i(t)}_t$. If $i(t)=0$, put
$$A^{i(t)}_t:=\{\rho<\kappa\mid  t^\smallfrown \langle \rho, H_t(\rho)\rangle^\smallfrown \langle A^*_\rho, H^*_\rho\rangle \parallel \varphi\}.$$
Otherwise, define
$$A^{i(t)}_t:=\{\rho<\kappa\mid \forall a\geq_{\mathbb{Q}(\rho)} H_t(\rho)\, (t^\smallfrown \langle \rho, a\rangle^\smallfrown \langle A^*_\rho, H^*_\rho\rangle \nparallel \varphi)\}.$$
Set $B^{i(t)}_t:=\{\rho<\kappa\mid \bar{H}(\rho)\geq_{\mathbb{Q}(\rho)} H_t(\rho)\}\cap A_t^{i(t)}.$ Clearly, $B^{i(t)}_t\in\mathcal{U}.$

Let $A':=\dom(\bar{H})\cap \diagonal_{t\preceq s}B^{i(t)}_t\in \mathcal{U}$ and set $q':=s^\smallfrown \langle A', H'\rangle,$ where $H':=\bar{H}\restriction A^{**}.$ We next show that $q'$ is almost the desired condition. 

\smallskip

Suppose that $t\preceq s$, $\rho\in A'_{\max(t)}$ and $a\geq_{\mathbb{Q}(\rho)} H'_{\max(t)}(\rho)$ are such that $t^\smallfrown \langle \rho, a\rangle^\smallfrown \langle A'_\rho, H'_{\rho}\rangle\parallel \varphi$. Since $\rho\in A'_{\max(t)}$ then $\rho\in B^{i(t)}_{t}\cap \dom(\bar{H})$, and so
$$a\geq_{\mathbb{Q}(\rho)} H'_{\max(t)}(\rho)=H'(\rho)=\bar{H}(\rho)\geq_{\mathbb{Q}(\rho)} H_t(\rho).$$
By the definition of $A^{i(t)}_t$  the above yields $i(t)=0$, hence $A'_{\max(t)}\s A^{0}_t$. Thereby,
$ t^\smallfrown \langle \rho, H_t(\rho)\rangle^\smallfrown \langle A^*_\rho, H^*_\rho\rangle \parallel \varphi$ for all $\rho\in A'_{\max(t)}.$
Finally, $(q'+t)\cat \rho\geq^{*} t^\smallfrown \langle \rho, H_t(\rho)\rangle^\smallfrown \langle A^*_\rho, H^*_\rho\rangle$, hence the former condition also decides $\varphi.$

\smallskip

Let us now $\leq^{**}$-extend $q'$ to $q^{*}$ to get a uniform decision about $\varphi.$

For each $t\preceq s$ put
$$A^{**}_{0,t}:=\{\rho\in A'_{\max(t)}\mid (q^*+t)\cat \rho\forces\varphi\},$$
and
$$A^{**}_{1,t}:=\{\rho\in A'_{\max(t)}\mid (q^*+t)\cat \rho\forces\neg\varphi\}.$$
Since this is a partition of $A'_{\max(t)}$ in two disjoint sets there is $i(t)\in \{0,1\}$ such that $A^{**}_{i(t),t}\in \mathcal{U}$. Put $A^{**}:=A'\cap \diagonal_{t\preceq s} A^{**}_{i(t),t}$ and $q^*:=s^\smallfrown \langle A^{**}, H^{**}\rangle$, where $H^{**}:=H'\restriction A^{**}.$ Arguing as above it is routine to check that $q^*$ has the desired property.
\end{proof}
Let $q^*:=s^\smallfrown \langle A^{**}, H^{**}\rangle$ be the condition provided by the above claim. Before proving the upcoming claim and completing the proof of the lemma let us note the following: if $u,v$ are conditions in $\mathbb{P}$ such that $u\parallel \varphi$ and $\stem(u)=\stem(v)$ then $v\parallel\varphi$, as well. In effect, if $w\geq v$ then there is a canonical way to produce a condition $z\geq w,u$ hence, in particular, $z\parallel \varphi.$

\begin{claim}
The condition $q^*$ decides $\varphi.$
\end{claim}
\begin{proof}
Suppose otherwise. Let $r\geq q^*$ be a condition with minimal $\ell(r)$ such that $r\parallel \varphi$. By the previous comments, $\ell(r)>\ell(q^*)$, for otherwise $q^*$ would decide $\varphi$.\footnote{Actually, it will decide $\varphi$ in the same way as $r$ does.} Let $t\preceq s$, $\rho\in A^{**}$ and $a\geq_{\mathbb{Q}(\rho)} H^{**}(\rho)$ be such that $r=t^\smallfrown \langle \rho, a\rangle^\smallfrown \langle B, L\rangle$. Once again, $t^\smallfrown \langle \rho, a\rangle^\smallfrown \langle A^{**}_\rho, H^{**}_{\rho}\rangle$ must decide $\varphi$, hence the above claim yields $(q^*+t)\cat \rho\parallel\varphi$ for all $\rho\in A^{**}_{\max(t)}$, and the corresponding decision is uniform. Since this collection of conditions forms a maximal antichain above $q^*+t$ we infer that  $q^*+t\parallel \varphi.$ From altogether we conclude that there is a condition $u$ (i.e., $q^*+t$) with $\ell(u)<\ell(r)$ which nevertheless decides $\varphi.$ This produces the desired contradiction.
\end{proof}
We have accomplished the proof of the lemma.
\end{proof}
The next simple lemma addresses the cardinal structure of $V^\mathbb{P}$. 
\begin{lemma}\label{LemmaCardinals} Let $G\s \mathbb{P}$ a generic filter and $\langle \rho_n\mid n<\omega\rangle$ be the induced Prikry sequence. Then the only cardinals in $V[G]$ are $V$-cardinals in the set $[\aleph_0, \rho_0]\cup[\kappa,\infty)$ and $$ \{\rho^{+k}_n\mid  n<\omega, \, k\in \{0,1,2,3\}\}.$$
By further forcing with $\col(\omega_1,{<}\rho_0)$ over $V[G]$  we have that $\rho_n:=\aleph_{4\cdot n+2}.$
Furthermore, the $\gch$ pattern is as follows:
\begin{itemize}
    \item $2^{\aleph_{4\cdot n+3}}=\aleph_{4\cdot n+5}$ for every $n<\omega$.
    \item For every cardinal $\lambda\notin\{\aleph_{4\cdot n+3}\mid n<\omega\}$, $2^\lambda=\lambda^+$.
\end{itemize}
\end{lemma}
\begin{proof}
Preservation of cardinals ${>}\kappa$ follows from the $\kappa^+$-cc of $\mathbb{P}$. Also, $\kappa$ will be preserved as a consequence of the preservation of the $\rho_n$'s. Cardinals ${\leq}\rho_0$ are preserved because the poset $\mathbb{P}/\langle \rho_0,a_0,A_0,H_0\rangle$ has the Prikry property and is $\rho_0^+$-closed with respect to $\leq^*$.\footnote{Here $\langle \rho_0,a_0,A_0,H_0\rangle$ is a condition in $G$.} Finally, let $n<\omega$ and $p\in G$ be with $\ell(p)=n+2$. By Lemma~\ref{factoring} we have that $\mathbb{P}/p$ is isomorphic to $\prod_{i\leq n}\mathbb{Q}(\rho_i,\rho_{i+1})\times \mathbb{P}/\langle \rho_{n+1}, a_{n+1}, A^p, H^p\rangle.$ The second of these factors preserves cardinals ${\leq}\rho_{n+1}$. Also, it does preserve $\gch_{\leq \rho_{n+1}}$.  Besides, since $\rho_{n}$ was inaccessible, $\prod_{i<n}\mathbb{Q}(\rho_i,\rho_{i+1})$ is a $\rho_n$-cc poset of size  $\rho_n$. In particular, this forcing preserves both cardinals ${\geq}\rho_n$ and  $\gch_{\geq\rho_n}$. Finally, in the resulting generic extension the poset  $\mathbb{Q}(\rho_n,\rho_{n+1})$ preserves $\rho_n^{+k}$ for $k\in\{0,1,2,3\}.$ This completes the verification and the $\mathbb{S}(\rho_n^+,\rho_{n}^{+3})$ makes $2^{\rho_n^+}=\rho_n^{+3}$.
\end{proof}

\begin{theorem}
Assume the $\gch$ holds and that there is a measurable cardinal $\kappa$. Then, there is a $\kappa^+$-cc generic extension where:
\begin{enumerate}
    \item $\aleph_\omega$ is strong limit;
    \item $\gch_{\geq \aleph_\omega}$ holds, hence $\Gal{\lambda^+}{\lambda^{++}}$ holds, for all $\lambda\geq \aleph_\omega$;
    \item  $\Gal{\aleph_{4\cdot(n+1)}}{\aleph_{4\cdot(n+1)+1}}$ fails for $n<\omega$. 
\end{enumerate}
\end{theorem}
\begin{proof}
Let us force with $\mathbb{P}\ast \lusim{\col}(\omega_1,{<}\lusim{\rho}_0)$. 
Let $G\s \mathbb{P}$ be a $V$-generic filter and $H\s \col(\omega_1,{<}\rho_0)$ be $V[G]$-generic. Clauses~(1) and (2) are evident. The argument for (3) is mutatis mutandi the same as that of Lemma~\ref{LemmaCardinals}. Note that it suffices to show that $\Gal{\rho_n^{++}}{\rho_{n+1}}$ fails in $V[G]$; indeed, $\col(\omega_1,{<}\rho_0)$ is $\rho_0$-cc  in $V[G]$ and thus has no effect upon this configuration. For details, see \cite[Lemma~2.1]{NegGalSing}.

Denote by $\langle \rho_n\mid n<\omega\rangle$ the Prikry sequence inferred from $G$. Let $n<\omega$ and $p\in G$ be such that $\ell(p)=n+2.$ By Lemma~\ref{factoring} we have that $\mathbb{P}/p$ is isomorphic to $\prod_{i\leq n}\mathbb{Q}(\rho_i,\rho_{i+1})\times \mathbb{P}/\langle \rho_{n+1}, a_{n+1}, A^p, H^p\rangle.$ Observe that the second of these factors has the Prikry property and is $\rho_{n+1}^{+}$-closed with respect to the $\leq^*$-ordering. Note that this is true even in the generic extension by the first factor in that this latter is a small forcing (i.e., of size $\rho_{n+1}$). In particular, the failure of $\Gal{\rho_n^{++}}{\rho_{n+1}}$ in $V[G]$ depends just on the effect of the first forcing. In this respect, since $\rho_{n}$ is an inaccessible cardinal, $\prod_{i<n}\mathbb{Q}(\rho_i,\rho_{i+1})$ has size $\rho_{n}$, so that $\gch_{\rho^+_n}$ will still hold in the resulting generic extension, $W$. Hence, over $W$, the poset  $\mathbb{Q}(\rho_n,\rho_{n+1})$ yields the failure of $\Gal{\rho^{++}_n}{\rho_{n+1}}.$ This completes the proof.
\end{proof}
Abraham and Shelah method produces the so-called \textit{ultimate failure}. Following \cite[\S2]{NegGalSing}, we say that Galvin's property \emph{ultimately fails} (say) at $\aleph_{n+1}$ if $\Gal{\aleph_{n+1}}{2^{\aleph_{n+1}}}$ fails. The next proposition suggests that combining Abraham and Shelah method with Prikry-type forcings seems unlikely to get infinitely-many consecutive failures of Galvin's property:
\begin{proposition}
If for some $1\leq n<\omega$, $\Gal{\aleph_{m}}{2^{\aleph_m}}$ fails for all $m\in [n,\omega)$ then $2^{\aleph_n}>\aleph_\omega$. In particular, $\aleph_\omega$ is not strong limit.
\end{proposition}
\begin{proof}
Assume otherwise. Let $n<\omega$ be such that $\Gal{\aleph_{m+1}}{\aleph_{m+2}}$ fails for all $m\in[n,\omega).$ Suppose towards a contradiction that $2^{\aleph_n}<\aleph_\omega$. Arguing as in \cite[Claim~2.15]{NegGalSing} one can show that $2^\theta=2^{\aleph_n}$ for all $\theta\in[\aleph_n, 2^{\aleph_n})$ (note that the $\gch$ must fail as Galvin's property does). Indeed, if $\theta$ is the minimal such that $2^\theta>2^{\aleph_n}$ then $\theta=\aleph_{m+1}$ and $2^{\aleph_m}=2^{\aleph_n}$. Hence the weak diamond $\Phi_{\aleph_{m+1}}$ holds, and thus $\Gal{\aleph_{m+1}}{2^{\aleph_{m+1}}}$ holds, as well (see \cite{MR3604115}). By our departing assumption this latter is certainly impossible. In particular, $2^{\aleph_n}=2^{\aleph_{n+1}}$. Similarly, we can now prove that for every $m\geq n$, $2^{\aleph_m}=2^{\aleph_{m+1}}$ and thus $2^{\aleph_n}\geq \aleph_\omega$, contradiction.
\end{proof}
\begin{question}\label{QuestionAboutManyfailures}
Is it consistent to have infinitely many consecutive failures of Galvin's property? If so, is the theory
$$\text{$``\mathsf{ZFC}+\forall n<\omega \;\neg \Gal{\aleph_{n+1}}{\aleph_{n+2}}$''}$$
consistent?
\end{question}
\section{On consistency strength}\label{SectionConsistencyStrength}
In this section we compute the exact consistency strength of 
 $$(\star) \ \  \ \ \   \text{``}\exists \kappa\, (\kappa\text{ strong limit, }
 \cf(\kappa)=\omega\;\&\;\Gal{\kappa^+}{\kappa^{++}}\text{ fails)''}$$
 and give a close-to-optimal upper bound for
 $$(\star\star) \ \ \ \ \ \ \ \ \ \  \text{``}\forall \kappa\, (\kappa\text{  singular}\;\Rightarrow\;\Gal{\kappa^+}{\kappa^{++}}\text{ fails)''.} \ \ \ \ \ \ \ \ \ \ \ $$  
The first answers \cite[Question~5.12]{NegGalSing} and the second improves the large-cardinal assumptions used in \cite[Theorem~2.3]{NegGalSing}.

Loosely speaking, the idea is to force a normal measure $\mathscr{U}$ on $\kappa$  concentrating on inaccessible cardinals $\alpha$ where $\Gal{\alpha^+}{\alpha^{++}}$ fails. As in \S\ref{MeasureConcentratingOnfailures}, this will be accomplished by forcing with an Easton-supported iteration of  $\mathbb{S}(\alpha,\alpha^{++})$'s at every inaccessible cardinal $\alpha\leq \kappa$. As the reader may have noticed, the difference now is that we also need to force at $\kappa$ and, as a result, the lifting arguments of \S\ref{MeasureConcentratingOnfailures} are not longer straightforward. The crux of the matter is lifting the relevant embedding  after forcing with $\mathbb{S}(\kappa,\kappa^{++})$. For this purpose, we  use an improvement of \emph{Woodin's surgery method} discovered by Ben-Shalom \cite{BenShalom} (see page~\pageref{BenShalom}). The resulting lifting can be shown to be the ultrapower by a normal measure over $\kappa$,  and  $\Gal{\kappa^+}{\kappa^{++}}$ fails. Finally,  one uses Prikry forcing with respect to this measure.

\smallskip

In Lemma~\ref{LemmaConstructingPair} we spell out the details of this construction, also dealing with the additional caveat of preserving $(\kappa+2)$-strongness of $\kappa$. This is the main ingredient to get a close-to-optimal bound for $(\star\star)$. Finally, in Theorem~\ref{exactstrengthlocal} we get the exact consistency strength for $(\star)$, even for $\aleph_{\omega}.$

\smallskip
For the rest of the section $\mathbb{P}$ will denote the Easton-support iteration forcing with $\mathbb{S}(\alpha,\alpha^{++})$ when $\alpha$ is an inaccessible cardinal ${\leq}\kappa$, for a given (large) cardinal $\kappa$. In the other stages we simply force with the trivial poset.

\subsection{Improving the consistency strength of the global failure}
The next lemma improves the large-cardinal assumptions of \cite[Lemma~2.4]{NegGalSing}:
\begin{lemma}\label{LemmaConstructingPair}
	Assume the $\mathsf{GCH}$ and that $\kappa$ is a $(\kappa+3)$-strong cardinal. Then, in $V^\mathbb{P}$ there is $(j,F)$ a weak constructing pair with $j$ witnessing that $\kappa$ is $(\kappa+2)$-strong.
\end{lemma}
\begin{proof}		
Work in $V$. Let $j\colon V\rightarrow M$ be an elementary embedding witnessing that $\kappa$ is  $(\kappa+3)$-strong. Without loss of generality assume that $j$ is the  ultrapower embedding by a $(\kappa,\kappa^{+3})$-extender $E$. Let $\ell\colon V\rightarrow \overline{M}$ be the ultrapower by the $(\kappa,\kappa^{++})$-extender induced by $E$. 
As usual, we denote by $k\colon \overline{M}\rightarrow M$  the factor map between $j$ and $\ell$; namely,  $k(\ell(f)(a)):=j(f)(a)$ for each $a\in[\kappa^{++}]^{<\omega}$. 
Note that $\crit(k)=\kappa^{+3}_{\overline{M}}$ and that $k$ has width $\kappa^{+3}_{\overline{M}}$.\footnote{
In addition,   $\text{$\kappa^{+i}_{\overline{M}}=\kappa^{+i}=\kappa^{+i}_M,\;$ for $i\in\{0,1\}$}.$}

Let $i\colon V\rightarrow N$ be the ultrapower embedding induced by $\ell$, and $\overline{k}$ the factor map between $\ell$ and $i$: i.e., $\bar{k}(i(f)(\kappa)):=j(f)(\kappa)$. Again, it is easy to verify that $\crit(\overline{k})=\kappa^{++}_N$ and that $\overline{k}$ has width $\kappa^{++}_N$.

\smallskip

Let $G\ast g$ be generic for $\mathbb{P}_\kappa\ast \lusim{\mathbb{S}}(\kappa,\kappa^{++})$ over $V$. Working over $V[G]$, observe that $\mathbb{S}(\kappa,\kappa^{++})_{N[G]}=\mathbb{S}(\kappa,\kappa^{++}_N)_{V[G]}$, as $N[G]$ is closed under $\kappa$-sequences in $V[G]$. Put $g_0:=g\cap \mathbb{S}(\kappa,\kappa_{N}^{++})_{V[G]}$. It is not difficult to check that this latter yields a generic for $\mathbb{S}(\kappa,\kappa_{N}^{++})_{V[G]}$ over $V[G]$.

By standard lifting arguments the embedding $i$ lifts to  another embedding $i\colon V[G]\rightarrow N[G\ast g_0\ast H]$, where $H\in V[G\ast g_0]$ is generic for the tail forcing $i(\mathbb{P})/(G\ast g_0)$. Similarly, $\overline{k}$ lifts to $\overline{k}\colon N[G\ast g_0\ast H]\rightarrow \overline{M}[G\ast g\ast \overline{k}``H]$: In effect, on one hand, $\crit(\overline{k})=\kappa^{++}_N$ and thus $\overline{k}``g_0\s g$; on the other hand, the width of $\overline{k}$ is $\kappa^{++}_N$ and $i(\mathbb{P})/(G\ast g_0)$ is $\kappa^{+3}_N$-closed in $N[G\ast g_0]$, hence $\overline{k}``H$ induces a generic filter for $\ell(\mathbb{P})/(G\ast g)$ over $\overline{M}[G\ast g]$. 
A similar argument shows that $k\colon \overline{M}\rightarrow M$ lifts to $k\colon \overline{M}[G\ast g\ast \overline{k}``H]\rightarrow M[G\ast g\ast (k\circ \overline{k})``H].$

All in all, in $V[G\ast g]$, we have a commutative diagram of embeddings given by $j\colon V[G]\rightarrow M[j(G)]$, $\ell\colon V[G]\rightarrow\overline{M}[\ell(G)]$, $i\colon V[G]\rightarrow N[i(G)]$, $k\colon \overline{M}[\ell(G)]\rightarrow M[j(G)]$ and $\overline{k}\colon N[i(G)]\rightarrow \overline{M}[\ell(G)]$.

\begin{claim}\label{ClaimGenericfilters}
$V[G\ast g]$ contains the following objects:
\begin{enumerate}
	\item $F\s \col(\kappa^{+4},{<}i(\kappa))_{N[i(G)]}$ a generic filter  over $N[i(G)]$;
	\item $K\s \Add(i(\kappa),i(\kappa))_{N[i(G)]}$ a generic filter over $N[i(G)]$.
\end{enumerate}

\end{claim}
\begin{proof}[Proof of claim]
We just prove (1) as (2) can be verified in the very same way. 
	Let us begin noticing that $i(\kappa)$ is an inaccessible  cardinal in $N[i(G)]$.  In particular, the poset $\col(\kappa^{+4},{<}i(\kappa))_{N[i(G)]}$ is $i(\kappa)$-cc and so there are at most $i(\kappa)^{<i(\kappa)}=i(\kappa)$-many maximal antichains.  Also, $$|i(\kappa)|=|\{i(f)(\kappa)\mid f\colon \kappa\rightarrow\kappa,\,f\in V[G]\}|,$$
	so that $|i(\kappa)|=|\kappa^\kappa|^{V[G]}=\kappa^+$. Finally,  standard diagonalization arguments yield the desired generic filter $F$ in $V[G\ast g]$.
\end{proof}
Since the width of $\overline{k}$ is $\kappa^{++}_N$, both $F$ and $K$ can be transferred to  generic filters for $\col(\kappa^{+4},{<}\ell(\kappa))_{\overline{M}[\ell(G)]}$ and $\Add(\ell(\kappa),\ell(\kappa))_{\overline{M}[\ell(G)]}$, respectively.  For simplicity, let us denote these generics also by $F$ and $K$.

\begin{claim}
$\ell$ lifts to $\ell\colon V[G\ast g]\rightarrow \overline{M}[\ell(G\ast g)]$ in $V[G\ast g].$	
\end{claim}
\begin{proof}
Work in $V[G]$. Factor $\mathbb{S}(\kappa,\kappa^{++})$ as $\Add(\kappa,\kappa^{++})\ast \lusim{\mathbb{Q}}$, where the latter stands for the quotient forcing. By arguments of Abraham and Shelah (see Lemma~\ref{Distributivityofthequotient}) the weakest condition of $\Add(\kappa,\kappa^{++})$ forces $\lusim{\mathbb{Q}}$ to be $\kappa^+$-distributive. For future convenience, let us split the generic filter $g$  as $g^c\ast g^q$.

We begin lifting $\ell$ via $\Add(\kappa,\kappa^{++})_{V[G]}$. For this, let $K$ be as before. Appealing to arguments of Ben Shalom \cite{BenShalom}\label{BenShalom} we can find a one-to-one  map $\varphi\colon \ell(\kappa)\rightarrow \ell(\kappa^{++})$ in $V[G\ast g]$ such that the function $K\diamond \varphi$ defined as  
$$K\diamond \varphi:=\{\langle\langle \varphi(\alpha),\beta\rangle, \gamma \rangle\mid \langle \langle \alpha,\beta\rangle, \gamma\rangle \in K\;\text{and}\;\alpha\in \dom(\varphi) \}$$
defines a generic for $\Add(\ell(\kappa),\ell(\kappa^{++}))_{\overline{M}[\ell(G)]}$. Evidently, $K\diamond \varphi\in V[G\ast g]$. 

 Using \emph{Woodin's surgery} (see \cite{Cummings-handbook}) one can alter (in $V[G\ast g]$) $K\diamond \varphi$ to a $\overline{M}[\ell(G)]$-generic filter $\overline{g}^c$ such that $\ell``g^c\s  \overline{g}^c$. Thus,  $\ell$ lifts to  $$\ell\colon V[G\ast g^c]\rightarrow \overline{M}[\ell(G)\ast \overline{g}^c].$$
 
 To complete the argument observe that $\ell$ is an embedding with width $\kappa$ and that $g^q$ is a generic for a $\kappa^+$-distributive forcing over $V[G\ast g^c]$. In particular, $\ell``g^q$ induces a generic filter for $\ell(\lusim{\mathbb{Q}})_{\overline{M}[\ell(G)\ast \overline{g}^c]}$. From altogether we conclude that $\ell$ lifts to $\ell\colon V[G\ast g]\rightarrow \overline{M}[\ell(G\ast g)]$.
\end{proof}
In addition, $k$ lifts to $k\colon \overline{M}[\ell(G\ast g)]\rightarrow M[j(G\ast g)]$, for $k$ has width $\kappa^{+3}_{\overline{M}}$ and $\mathbb{S}(\ell(\kappa),\ell(\kappa^{++}))_{\overline{M}[\ell(G)]}$ is $\ell(\kappa)$-closed in ${\overline{M}[\ell(G)]}$. Incidentally, by commutativity of the diagram, $j$ lifts to $j\colon V[G\ast g]\rightarrow M[j(G\ast g)]$.

\smallskip

Write $V^*:=V[G\ast g]$, $M^*:=M[j(G\ast g)]$ and $\overline{M}^*:=\overline{M}[\ell(G\ast g)]$.
\begin{claim}\label{Claimfinishingconstructing}\hfill
	\begin{enumerate}
		\item $j\colon V^*\rightarrow M^*$ witnesses that $\kappa$ is $(\kappa+2)$-strong;
		\item $\ell\colon V^*\rightarrow \overline{M}^*$ is the ultrapower embedding derived from $j$;
		\item $F$ is $\overline{M}^*$-generic for $\col(\kappa^{+4},{<}\ell(\kappa))_{\overline{M}^*}$.
	\end{enumerate}
	In particular, $(j,F)$ is a weak constructing pair.
\end{claim}
\begin{proof}[Proof of claim]
(1) It suffices to check the following two bullets: 

$\br$ \underline{$M^*$ is closed under $\kappa$-sequences in $V^*$:}  First, $\mathbb{P}_\kappa\ast \lusim{\Add}(\kappa,\kappa^{++})$ is $\kappa^{+}$-cc and so $M[G\ast g^c]$ is closed under $\kappa$-sequences in $V[G\ast g^c]$. Second, the quotient forcing $\mathbb{Q}$ is $\kappa^{+}$-distributive in $V[G\ast g^c]$ and so $M[G\ast g]$ is $\kappa$-closed in $V[G\ast g]$. Finally, the tail forcing $\ell(\mathbb{P})/(G\ast g)$ is $\kappa^{+}$-closed in $V[G\ast g]$ and thus $M^*=M[\ell(G\ast g)]$ is closed under $\kappa$-sequences in $V[G\ast g]$.

$\br$ \underline{$V^*_{\kappa+2}\s M^*$:} Let $x\in V^*_{\kappa+2}$. Since $\mathbb{P}$ is a $\kappa^{++}$-cc forcing notion then there is $\tau\in V_{\kappa+2}\s M$ such that $x=\tau_{G\ast g}$. Hence $x\in M^*$, as wanted. 
\smallskip

(2) Let $\iota\colon V^*\rightarrow\mathcal{N}$ be the ultrapower embedding inferred from $\ell$ and $\chi\colon \mathcal{N}\rightarrow \overline{M}^*$ be the corresponding factor map. Since $\crit(k\circ \chi)>\kappa$,  $\mathcal{N}$ is the ultrapower derived from $j$.  Thus, it suffices to argue that $\chi=\id$, which amounts to show that $\chi$ is surjective.\footnote{Recall that the identity is the unique isomorphism between transitive models of $\mathsf{ZFC}$.} Before checking this observe that $\crit(\chi)>\kappa^{++}$: In effect, $\kappa^+_\mathcal{N}=\kappa^+$ and also, since $\mathcal{N}\models ``2^\kappa=\kappa^{++}$'' and $\mathcal{P}(\kappa)\s \mathcal{N}$, $\kappa^{++}_\mathcal{N}=(2^\kappa)_\mathcal{N}\geq 2^\kappa=\kappa^{++}.$ Altogether, $\crit(\chi)>\kappa^{++}$.

 Let $x\in \overline{M}^*$. Since $\overline{M}$ was the ultrapower by a $(\kappa,\kappa^{++})$-extender it follows from standard forcing arguments (see e.g., \cite[Proposition~9.4]{Cummings-handbook}) that $$\overline{M}^*=\{\ell(f)(a)\mid f\in V[G\ast g],\; f\colon [\kappa]^{|a|}\rightarrow V[G\ast g],\; a\in [\kappa^{++}]^{<\omega}\}.$$
Thus, $x$ is of the form $\ell(f)(a)$ for some $f$ and $a$. Since $\crit(\chi)>\kappa^{++}$, $\kappa^{++}\s \ran(\chi)$ and hence there is $\bar{a}\in\mathcal{N}$ such that $\chi(\bar{a})=a$. All in all, $x=\ell(f)(a)=\chi(\iota(f))(\chi(\bar{a}))=\chi(\iota(f)(\bar{a}))$ and thus $x\in\ran(\chi)$, as wanted.

(3) Since $\ell(\mathbb{P})/\ell(G)$ is $\ell(\kappa)$-closed in $\overline{M}[\ell(G)]$ the poset $\col(\kappa^{+4},{<}\ell(\kappa))$ is computed both by $\overline{M}^*$ and $\overline{M}[\ell(G)]$ in the same manner. Besides, this is a $\ell(\kappa)$-cc forcing  hence any maximal antichain in $\overline{M}^*$ belongs to $\overline{M}[\ell(G)]$. Thereby, $F$ is $\overline{M}^*$-generic for the poset $\col(\kappa^{+4},{<}\ell(\kappa))_{\overline{M}^*}$
\end{proof}
This completes the proof of the theorem.
\end{proof}
The following is an immediate consequence of the above lemma:
\begin{theorem}
	Assume that $\mathsf{ZFC}$ is consistent with the existence of a $(\kappa+3)$-strong cardinal.
Then   $\mathsf{ZFC}$ is also  consistent with  
$$\text{$``{\rm Gal}(\mathrm{Cub}_{\kappa^+},\kappa^+,\kappa^{++})$ fails at every limit cardinal $\kappa$''.}$$
Moreover, $\mathsf{ZFC}$ is consistent with 
$$\text{$``\Gal{\aleph_{4\cdot\xi+1}}{\aleph_{4\cdot\xi+2}}$ fails for every $\xi\in \mathrm{Ord}$''.}$$
\end{theorem}
\begin{proof}
	Let $\kappa$ be a $(\kappa+3)$-strong cardinal. Appealing to Lemma~\ref{LemmaConstructingPair} we get a model of \textsf{ZFC} where $\Gal{\kappa^+}{\kappa^{++}}$ fails and there is a weak constructing pair $(j,F)$ witnessing that $\kappa$ is $(\kappa+2)$-strong. Denote this model by $V$ and let $u_*$ be the sequence inferred from $(j,F)$. Arguing  as in \cite[Lemma~1]{CumGCH}, for each 
$\alpha<\kappa^{+3}$ the sequence $u_*\restriction\alpha$ exists and belongs to $\mathcal{U}_\infty$.  
In particular, there is $\alpha<\kappa^{+3}$ such that the sequence $u:=u_*\restriction\alpha\in \mathcal{U}_\infty$ has a repeat point. From this point on argue exactly as in \cite[Theorem~2.3]{NegGalSing}.
\end{proof}

\subsection{The exact consistency strength for the local failure}
In this section we compute the exact consistency strength of the failure of Galvin's property at the successor of a strong limit singular cardinal.
\begin{theorem}\label{exactstrengthlocal}
Assume that there is a cardinal $\kappa$  carrying a $(\kappa,\kappa^{++})$-extender. Then, there is a  generic extension of the set-theoretic universe where  $\aleph_\omega$ is strong limit and $\Gal{\aleph_{\omega+1}}{\aleph_{\omega+2}}$ fails.
\end{theorem}
\begin{proof}
Let $j\colon V\rightarrow M$ be the ultrapower induced by a $(\kappa,\kappa^{++})$-extender, and $i\colon V\rightarrow N$ be the natural ultrapower embedding inferred from $j$. Arguing as in  Lemma~\ref{LemmaConstructingPair} lift $j$ to $j\colon V[G\ast g]\rightarrow M[j(G\ast g)]$. In addition, find $F\s \col(\kappa^{+3},{<}j(\kappa))_{M[j(G)]}$ 
a generic filter over $M[j(G)]$ living in $V[G\ast g].$ Since $\mathbb{S}(j(\kappa),j(\kappa^{++}))_{M[j(G)]}$ is $j(\kappa)$-closed we can argue as in Claim~\ref{Claimfinishingconstructing}(3) that $F$ is generic for $\col(\kappa^{+3},{<}j(\kappa))_{M[j(G\ast g)]}$ over $M[j(G\ast g)].$
By  Claim~\ref{Claimfinishingconstructing}(2), $M[j(G\ast g)]$ is just the ultrapower by the normal measure inferred from $j$.

Working in $V[G\ast g]$, let $\mathcal{U}$ denote this  measure and $\mathbb{Q}$ be the Prikry forcing with interleaved collapses relative to $\mathcal{U}$ and the \emph{guiding generic} $F$. Since in $V[G\ast g]$ the principle $\Gal{\kappa^+}{\kappa^{++}}$ fails and $\mathbb{Q}$ is $\kappa^+$-cc, it also fails in any generic extension of $V[G\ast g]$ by $\mathbb{Q}.$ Finally, observe that in this latter model $\aleph_\omega$  is strong limit and $\kappa=\aleph_\omega$.
\end{proof} 

\begin{corollary}
 $``\mathsf{ZFC} +  \text{$\aleph_\omega$ is strong limit  $+ \Gal{\aleph_{\omega+1}}{\aleph_{\omega+2}}$ fails}$'' and $``\mathsf{ZFC} + \exists \kappa\,(o(\kappa)=\kappa^{++})$'' are equiconsistent theories.
\end{corollary}

 \bibliographystyle{alpha}
\bibliography{arlist}
\end{document}